\documentclass[twoside]{article} 

\usepackage{jmlr2e} 
\usepackage{algorithm}
\usepackage{algorithmic}
\usepackage{hyperref}
\usepackage{natbib}
\usepackage{lipsum}
\usepackage{setspace}

\usepackage{amsmath,amssymb}
\usepackage{pgfplots}
\usepackage{color}
\usepackage{flushend}
\usepackage{multirow}
\usepackage{rotating}
\usepackage{appendix}
\usepackage{tikz,pgfplots}
\usepackage{subcaption}
\usepackage{bbm}

\input{mysymbol.sty}
\usepackage{needspace}

\newtheorem{assumption}{\hspace{0pt}\bf Assumption}

\newcommand{\QED}{$\blacksquare$}

\newcommand{\graph}{\mathcal{G}}
\newcommand{\nodes}{\mathcal{N}}
\newcommand{\edges}{\mathcal{E}}

\begin{document}

% Heading arguments are {volume}{year}{pages}{submitted}{published}{author-full-names}

\jmlrheading{1}{2000}{1-48}{4/00}{10/00}{Aryan Mokhtari, Hamed Hassani, and Amin Karbasi}

% Short headings should be running head and authors last names

\ShortHeadings{Decentralized Submodular Maximization}{Mokhtari, Hassani, and Karbasi}
\firstpageno{1}

\title{Decentralized Submodular Maximization:\\ Bridging Discrete and Continuous Settings}

\name{Aryan Mokhtari$^{\S}$, Hamed Hassani$^{\dag}$, and Amin Karbasi$^{ \ddagger}$ \thanks{This work was done while A. Mokhtari was a Research Fellow at the Simons Institute for the Theory of Computing.}
\address{\normalsize $^\S$Laboratory for Information and Decision Systems, Massachusetts Institute of Technology \\
\normalsize$^{\dag}$Department of Electrical and Systems Engineering, University of Pennsylvania\\
\normalsize$^{\ddagger}$Department of Electrical Engineering and Computer Science, Yale University\\
\normalsize aryanm@mit.edu, hassani@seas.upenn.edu, amin.karbasi@yale.edu}}

\maketitle
%\thispagestyle{empty}
% this must go after the closing bracket ] following \twocolumn[ ...

% This command actually creates the footnote in the first column
% listing the affiliations and the copyright notice.
% The command takes one argument, which is text to display at the start of the footnote.
% The \icmlEqualContribution command is standard text for equal contribution.
% Remove it (just {}) if you do not need this facility.

%\printAffiliationsAndNotice{}  % leave blank if no need to mention equal contribution
%\printAffiliationsAndNotice{\icmlEqualContribution} % otherwise use the standard text.

%!TEX root = example_paper.tex
\begin{abstract}
In this paper, we showcase the interplay between discrete and continuous optimization in network-structured settings. We propose the first  fully decentralized optimization method for a wide class of non-convex objective functions that possess a diminishing returns property. More specifically, given an arbitrary connected network and a global \textit{continuous} submodular function, formed by a sum of local functions, we develop \alg (DCG), a message passing algorithm that converges to the tight $(1-1/e)$ approximation factor of the optimum global solution using only local computation and communication. We also provide strong convergence bounds as a function of network size and spectral characteristics of the underlying topology. Interestingly,
DCG readily provides a simple recipe for decentralized  \textit{discrete} submodular maximization through the means of continuous relaxations. Formally, we demonstrate that by lifting the local discrete functions to continuous domains and using DCG as an interface we can develop a consensus algorithm that also achieves the tight $(1-1/e)$ approximation guarantee of the global discrete solution once a proper rounding scheme is applied. 

%
%In this paper, we propose the first fully decentralized optimization method for a wide class of non-convex objective functions that diminishing returns property. Given an arbitrary network and a global \textit{continuous} submodular function, formed by a sum of local functions, we develop XXX, a message passing algorithm that converges to the tight $(1-1/e)$ fraction of the the optimum solution using only local computation and communication. We provide strong convergence bounds as a function of network size and spectral characteristics of the underlying topology. More specifically, we show that the convergence rate  scales inversely proportional  in the spectral gap of the network. Interestingly,
%XXX readily provides a simple recipe for decentralized  \textit{discrete} submodular maximization. Even though the global objective function  is discrete, we demonstrate that by lifting the local discrete functions to continuous domains, we can develop a consensus algorithm when a proper rounding scheme is applied. Therefore, we showcase the interplay between discrete and continuous submodular maximization in network-structured optimization settings.
\end{abstract}
\newpage
%!TEX root = example_paper.tex
\section{Introduction}
In recent years, we have reached  unprecedented data volumes that are high dimensional and  sit over (clouds of) networked machines. As a result,  \textit{decentralized} collection of these data sets along with accompanying distributed optimization methods are not only desirable  but very often necessary \citep{BoydEtalADMM11}. 

The focus of this paper is on decentralized optimization, the goal of which is to maximize/minimize a global objective function  --distributed over a network of computing units-- through local computation and communications among nodes. A canonical example in machine learning is fitting models using  M-estimators where given a set of data points  the parameters of the model are estimated through an empirical risk minimization \citep{DBLP:books/daglib/0097035}.  Here, the global objective function is defined as an average of local loss functions associated with each data point. Such local loss functions can be convex (e.g., logistic regression, SVM, etc) or non-convex (e.g., non-linear square loss, robust regression, mixture of Gaussians, deep neural nets, etc) \citep{mei2016landscape}. Due to the sheer volume of data points,  these optimization tasks cannot be fulfilled on a single computing cluster node. Instead, we need to opt for decentralized  solutions that can efficiently exploit dispersed (and often distant) computational resources linked through a tightly connected network. Furthermore, local computations should  be light so that they can be done on single machines. In particular, when the data is high dimensional, extra care should be given to any optimization procedure that relies on projections over the feasibility domain.

In addition to large scale machine learning applications, decentralized optimization is a  method of choice in many other domains such as Internet of Things (IoT) \citep{abu2013data}, remote sensing \citep{ma2015remote},  multi-robot systems \citep{tanner2005towards}, and sensor networks \citep{rabbat2004distributed}. In such scenarios, individual  entities can communicate   over a network and  interact with the environment by exchanging the data generated through sensing. At the same time they can react to events and trigger actions to control the physical world. These applications highlight another important aspect of decentralized optimization where \textit{private} data is collected by different sensing units \citep{yang2017survey}. Here again, we aim to optimize a global objective function while avoiding to share the private data among computing units. Thus, by design, one cannot solve such private optimization problems in a centralized manner and should rely on decentralized solutions where local private computation is done where the data is collected. 

Continuous submodular functions, a broad subclass of non-convex functions with diminishing returns property, have recently received considerable attention \citep{bach2015submodular, bian16guaranteed}. Due to their interesting structures that allow strong approximation guarantees \citep{mokhtari2017conditional, bian16guaranteed}, they have found various applications, including robust budget allocation \citep{staib2017robust,soma2014optimal},  online resource allocation \citep{eghbali2016designing},  learning assignments \citep{golovin2014online}, as well as Adwords for e-commerce and advertising \citep{devanur2012online, mehta2007adwords}. However, all the existing work suffer from centralized computing. Given that many information gathering, data summarization, and non-parametric learning problems are inherently related to large-scale submodular maximization, the demand for a fully decentralized solution is immediate. In this paper, we develop the first decentralized framework for both continuous and discrete submodular functions.  Our contributions are as follows:
\begin{itemize}

\item \textit{Continuous submodular maximization:} For any  global objective function that is monotone and continuous DR-submodular and subject to any down-closed and bounded convex body, we develop \alg, a decentralized and \textit{projection-free} algorithm that achieves the tight $(1 - 1/e- \epsilon)$ approximation guarantee in $O(1/\epsilon^2)$ rounds of local communication. 
\vspace{-.1cm}
\item \textit{Discrete submodular maximization:} For any  global objective function that is monotone and submodular and subject to any matroid constraint, we develop a discrete variant of the \alg algorithm that achieves the tight $(1-1/e-\epsilon)$ approximation ratio in $O(1/\epsilon^3)$ rounds of communication. 
\end{itemize}

\section{Related Work}

%%%%%%%%%%%%%%%%%%%%%%%%%%%%%%%%%%%%%
%%%%%%%%%%%%%%%%%%%%%%%%%%%%%%%%%%%%%
%%%%%%      S  E  C  T  I  O  N     %%%%%%%%%%%%%%%%%%%
%%%%%%%%%%%%%%%%%%%%%%%%%%%%%%%%%%%%%
%%%%%%%%%%%%%%%%%%%%%%%%%%%%%%%%%%%%%
%%%%%%%%%%%%%%%%%%%%%%%%%%%%%%%%%%%%%

Decentralized optimization  is a challenging problem  as nodes only have access to separate components of the global objective function, while they  aim to collectively reach the global optimum point. Indeed, one naive approach to tackle this problem is to broadcast local objective functions to all the nodes in the network and then solve the problem locally. However, this scheme requires high communication overhead and disregards the privacy associated with the data of each node.  An alternative approach is the master-slave setting \citep{bekkerman2011scaling,shamir2014communication,zhang2015disco} where at each iteration, nodes use their local data   to compute the information needed by the master node. Once the master node receives all the local information, it updates its decision and broadcasts the decision to all the nodes. Although this scheme protects the privacy of nodes it is  not robust to machine failures and is prone to high overall communication time. In decentralized methods, these issues are overcame by removing the master node and considering each node as an independent unit that is allowed to exchange information with its neighbors.
%
%and avoids high communication cost of exchanging local functions, it is not robust to machine failures. In particular, failure of the master node completely breaks down the system. Furthermore, it could be inefficient in terms of overall communication time. For example, on a line-graph of size $n$ with link delays $\tau$, the communication time  required at each round to send the information of all the nodes to the master node is $\mathcal{O}(\tau n)$.  
%
%
%In this setting, and at each iteration, nodes use their local information (e.g., their local functions) to compute the information needed by the master node (e.g., their local gradients). Then, they send a low communication cost message (e.g., local variables or gradients) to the master node and it uses the received information to locally update its decision variable. At the end of each iteration, master broadcasts the updated decision variable to all the workers. Although this scheme protects the privacy of nodes and avoids high communication cost of exchanging local functions, it could be inefficient in terms of overall communication time when the graph is sparse. As an example, when the graph is a line, if the communication time for sending a particular information, e.g., local gradient, to a neighbor is $\tau$, then the communication time required at each round to send the information of all the nodes in the network to the master node is of $\mathcal{O}(\tau n)$ which grows linearly by the size of the network $n$. 
\vspace{.1cm}

{Convex decentralized consensus optimization} is a relatively mature area with a myriad of primal and dual algorithms \citep{bertsekas1989parallel}. Among primal methods, decentralized (sub)gradient descent  is perhaps the most well known algorithm which is a mix of local gradient descent and successive averaging \citep{nedic2009distributed,yuan2016convergence}. It also can be interpreted as a penalty method that encourages agreement among neighboring nodes. This latter
interpretation has been exploited to solve the penalized objective function using accelerated gradient descent \citep{jakovetic2014fast,qu2017accelerated}, Newton's method \citep{mokhtari2017network,bajovic2017newton}, or quasi-Newton algorithms \citep{eisen2017decentralized}. The methods that operate in the dual domain consider a constraint that enforces equality between nodes' variables and solve the problem by ascending on the dual function to find optimal Lagrange multipliers. A short list of dual methods are the alternating directions method of multipliers (ADMM) \citep{DBLP:journals/tsp/SchizasRG08,BoydEtalADMM11}, dual ascent algorithm \citep{rabbat2005generalized}, and augmented Lagrangian methods \citep{jakovetic2015linear,chatzipanagiotis2015convergence,mokhtari2016dsa}. Recently, there have been many attempts to extend the tools in decentralized consensus optimization to the case that the objective function is non-convex \citep{di2016next,sun2016distributed,hajinezhad2016nestt,tatarenko2017non}. However, such works are mainly concerned with reaching a stationary point and naturally cannot provide any optimality guarantee. 

\vspace{.1cm}

In this paper, our focus is to provide the first decentralized algorithms for both discrete and continuous submodular functions. Indeed, it is known that the centralized  greedy approach  of \citep{nemhauser1978analysis}, and its many variants  \citep{feige2011maximizing, buchbinder2015tight, buchbinder2014submodular, feldman2017greed,DBLP:conf/icml/MirzasoleimanBK16}, reach tight approximation guarantees in various scenarios.  Since such algorithms are sequential in nature, they do not scale to massive datasets. To partially resolve this issue, MapReduce style algorithms, with a master-slave architecture,  have been proposed \citep{DBLP:conf/nips/MirzasoleimanKSK13,kumar2015fast,DBLP:conf/icml/BarbosaENW15,DBLP:conf/stoc/MirrokniZ15, qu2015distributed}. 

%
%achieves the tight $(1-1/e)$ approximation guarantee for a monotone submodular set function subject to a cardinality constraint. To achieve the same approximation guarantee under a matroid constraint, it was observed that one needs to consider the multi-linear relaxation of the discrete submodular function and run a variant of the Frank-Wolfe algorithm, known as continuous greedy \citep{calinescu2011maximizing}. Recently,  greedy-like algorithms have been proposed to achieve tight approximation guarantees  for non-monotone and more general constraints \citep{feige2011maximizing, buchbinder2015tight, buchbinder2014submodular, feldman2017greed,DBLP:conf/icml/MirzasoleimanBK16}. Since such algorithms are sequential in nature, they do not scale to massive datasets. To partially solve this issue, MapReduce style algorithms, with a master-slave architecture,  have gained a lot of interest \citep{DBLP:conf/nips/MirzasoleimanKSK13,kumar2015fast,DBLP:conf/icml/BarbosaENW15,DBLP:conf/stoc/MirrokniZ15}. 

\vspace{.1cm}

One can  extend the notion of diminishing returns to continuous domains \citep{wolsey1982analysis,bach2015submodular}.  
%
%Such functions  arise   in many learning applications such as robust budget allocation \citep{staib2017robust,soma2014optimal},  online resource allocation \citep{eghbali2016designing},  learning assignments \citep{golovin2014online}, as well as Adwords for e-commerce and advertising \citep{devanur2012online, mehta2007adwords}. 
Even though continuous submodular functions are not generally convex (nor concave) \citet{hassani2017gradient} showed that in the monotone setting and subject to a general bounded convex body constraint, stochastic gradient methods can achieve a $1/2$ approximation guarantee. The approximation guarantee can be tightened to $(1-1/e)$ by using Frank-Wolfe \citep{bian16guaranteed} or stochastic Frank-Wolfe \citep{mokhtari2017conditional}. 
\vspace{.2cm}
\section{Notation and Background}
In this section, we  review the notation that we use throughout the paper. We then give the precise definition of submodularity in discrete and continuous domains.\\
\\
\textbf{Notation.} 
Lowercase boldface $\bbv$ denotes a vector and uppercase boldface $\bbW$ a matrix. The $i$-th element of  $\bbv$ is written as $v_i$ and the element on the $i$-th row and $j$-th column of  $\bbW$ is denoted by $w_{i,j}$. We use $\|\bbv\|$ to denote the Euclidean norm of vector $\bbv$ and $\|\bbW\|$ to denote the spectral norm of matrix $\bbW$. The null space of matrix $\bbW$ is denoted by $\rm{null}(\bbW)$. The inner product of vectors $\bbx,\bby$ is indicated by $\langle \bbx,\bby\rangle$, and the transpose of a vector $\bbv$ or matrix $\bbW$ are denoted by $\bbv^\dag$ and $\bbW^\dag$, respectively. The vector $\bbone_n\in \reals^n$ is the vector of all ones with $n$ components, and the vector $\bb0_p\in \reals^p$ is the vector of all zeros with $p$ components. \\
\\
\textbf{Submodulary.} A \textit{set} function $f:2^V\rightarrow \reals_+$, defined on the ground set $V$,  is called submodular if for all  $A,B\subseteq V$, we have $$f(A)+f(B)\geq f(A\cap B) + f(A\cup B).$$
We often need to maximize submodular functions subject to a down-closed set family $\mathcal{I}$. In particular, we say $\mathcal{I}\subset 2^V$ is a matroid if 1) for any $A\subset B\subset V$, if $B\in \mathcal{I}$, then $A\in \mathcal{I}$ and 2) for any $A,B\in \mathcal{I}$ if $|A|<|B|$, then there is an element $e\in B$ such that $A\cup\{e\}\in \mathcal{I}$.

The notion of submodularity goes beyond the discrete domain \citep{wolsey1982analysis, vondrak2007submodularity, bach2015submodular}.
Consider a continuous function $F: \ccalX \to \reals_{+}$ where the set $\ccalX \subseteq \mathbb{R}^p$ is of  the form $\ccalX=\prod_{i=1}^p\ccalX_i$ and each $\ccalX_i$ is a compact subset of $\reals_+$. We call the \textit{continuous} function $F$ submodular if for all $\bbx,\bby\in \ccalX$ we have
%%%
\begin{align}\label{eq:submodular_def}
F(\bbx) + F(\bby) \geq F(\bbx \vee	 \bby) + F(\bbx \wedge \bby) ,
\end{align}
%%%
where $\bbx \vee \bby := \max (\bbx ,\bby )$ (component-wise) and $\bbx \wedge \bby := \min (\bbx ,\bby )$ (component-wise). %It is not hard to verify that submodular functions are neither convex nor concave.  
In this paper, our focus is on differentiable  continuous submodular functions with two additional  properties:   monotonicity and diminishing returns.  Formally, a submodular function $F$ is monotone if 
%%%
\begin{align}\label{eq:monotone_def}
\bbx \leq \bby  \quad \Longrightarrow \quad F(\bbx) \leq  F(\bby),
\end{align}
%%%
for all $\bbx,\bby\in \ccalX$. Note that $\bbx \leq \bby$  in \eqref{eq:monotone_def} means that $x_i\leq y_i$ for all $i=1,\dots,p$. Furthermore, a differentiable submodular function $F$ is called \textit{DR-submodular} (i.e., shows diminishing returns) if the gradients are antitone, namely, for all $\bbx,\bby\in \ccalX$ we have 
%%%
\begin{align}\label{eq:antitone_def}
\bbx \leq \bby  \quad \Longrightarrow \quad \nabla F(\bbx) \geq \nabla  F(\bby).
\end{align}
When the function $F$ is twice differentiable, submodularity implies that all cross-second-derivatives are non-positive \citep{bach2015submodular},
%%%
%\begin{equation}
%\forall\ i\neq j,\ \ \forall\ \bbx\in \ccalX, ~~ \frac{\partial^2 F(\bbx)}{\partial x_i \partial x_j} \leq 0, 
%\end{equation}
%%%
 and DR-submodularity implies that all second-derivatives  are  non-positive \citep{bian16guaranteed}
%%%
%\begin{equation}
%\forall\ i,j,\ \ \forall\ \bbx\in \ccalX, ~~ \frac{\partial^2 F(\bbx)}{\partial x_i \partial x_j} \leq 0.
%\end{equation}
In this work, we consider the maximization of continuous submodular functions subject to \textit{down-closed convex bodies} $\ccalC \subset \reals_{+}^p$ defined as follows. For any two vectors $\bbx, \bby\in \reals_{+}^p$, where $\bbx\leq \bby$, down-closedness means that if $\bby\in \ccalC$, then so is $\bbx\in \ccalC$. Note that for a down-closed set we have $\mathbf{0}_p \in \ccalC$.
%
%Let us now explain the problem settings considered in this paper. 
%We provide in the next section our algorithms for these problems and discuss the
%assumptions we use in our convergence analysis. 
\section{Decentralized Submodular Maximization}
In this section, we  state the problem of decentralized submodular maximization in  continuous and discrete settings.\\ 
\\
\noindent\textbf{Continuous Case.} We consider a set of $n$ computing machines/sensors that communicate over a graph to maximize a global objective function.
Each machine can be viewed as a node $i \in \nodes \triangleq \{1,\cdots,n\}$. We further assume that the possible
communication links among nodes are given by a bidirectional connected \emph{communication graph} $\graph = (\nodes,\edges)$ where each node can only communicate with its neighbors in $\graph$. We formally use $\ccalN_i$ to denote node $i$'s neighbors.
%
%
%the set of nodes that node $i$ can communicate with, i.e., node $i$'s neighbors.
In our setting, we assume that each node $i \in \nodes$ has access to a local function $F_i:\ccalX\to\reals_{+}$. The nodes cooperate in order to maximize the aggregate monotone and continuous DR-submodular function $F:\ccalX\to\reals_{+}$ subject to a down-closed convex body $\ccalC\subset \ccalX\subset \reals_{+}^p$, i.e., 
%
% defined as
%\begin{equation}
%F(\bbx) := \frac{1}{n} \sum_{i=1}^{n}F_i(\bbx).
%\end{equation}
%More precisely, agents cooperate in solving the global optimization problem
\begin{equation}\label{original_optimization_problem1}
 \ \max_{\bbx\in \ccalC}  F(\bbx) 
       \  =\ \max_{\bbx\in \ccalC} \frac{1}{n}\sum_{i=1}^{n} F_i(\bbx).
\end{equation}
%where $\ccalC\subset \ccalX\subset \reals_{+}^p$ is a convex body. 
%As each local function $F_i$ is only accessible to $i$-th agent (node), 
%in order to maximize \eqref{original_optimization_problem1} the agents have to 
%communicate information about their local functions (possibly combined with other information) 
%to their neighbors, and such information will be propagated across the network  through many rounds of communication. 
The goal is to design a message passing algorithm to solve \eqref{original_optimization_problem1} such that:  (i) at each iteration $t$, the nodes send their messages (and share their information) to their neighbors in $\graph$, and (ii) as $t$ grows, all the nodes reach to a point $\bbx \in \ccalC$ that provides a (near-) optimal solution for  \eqref{original_optimization_problem1}.  \\
\\
\noindent \textbf{Discrete Case.}
Let us now consider the discrete counterpart of problem~\eqref{original_optimization_problem1}. In this setting, each node $i \in \nodes$ has access to a local  \emph{set} function
$f_i:2^V \to \mathbb{R}_+$.
The nodes cooperate in maximizing the aggregate monotone submodular  function $f:2^V \to\reals_{+}$ subject to a matroid constraint $\ccalI$, i.e.,
\begin{equation}\label{original_optimization_problem2}
\ \max_{S \in \mathcal{I}}  f(S) 
       \  =\ \max_{S\in \mathcal{I}} \frac{1}{n}\sum_{i=1}^{n} f_i(S).
\end{equation}
Note that even in the centralized case, and  under reasonable complexity-theoretic assumptions, the best approximation guarantee we can achieve for Problems~\eqref{original_optimization_problem1} and \eqref{original_optimization_problem2} is $(1-1/e)$ \citep{feige1998threshold}. In the following, we show that it is possible to achieve the same approximation guarantee in a decentralized setting. 
%defined as
%\begin{equation}
%f(S) := \frac{1}{n} \sum_{i=1}^{n}f_i(S),
%\end{equation}
%for $S \subseteq V$.
%More precisely, agents cooperate in solving the global optimization problem
%\begin{equation}\label{original_optimization_problem2}
%\text{OPT}\ = \ \max_{S \in \mathcal{I}}  f(S) 
%       \  =\ \max_{S\in \mathcal{I}} \frac{1}{n}\sum_{i=1}^{n} f_i(S),
%\end{equation}
%where $\mathcal{I}$ is a general matroid constraint (e.g. a cardinality constraint). 

%We consider the following problem of distributed constrained optimization:

%Let us denote the
%local decision variable associated with node $i \in \nodes$ as $x_i$. Each local
%decision variable is restricted to lie in a local constraint set
%.
%\begin{itemize}
%\item Submodular set functions and DR submodular set functions
%\item The decenteralized setting: the communication graph, node's abilities, precise formulation of the optimization problem
%\item Might be easier to make the connection to the discrete setting here but refer to a "Bridging..." section that will appear later on in the paper  
%\end{itemize}

%!TEX root = example_paper.tex
\section{Decentralized Continuous Greedy Method}

In this section, we introduce the \alg (DCG) algorithm for solving Problem \eqref{original_optimization_problem1}.
Recall that in a decentralized setting, the nodes have to cooperate (i.e., send messages to their neighbors) in order to solve the global optimization problem. We will explain how such messages are designed and communicated in DCG. Each node $i$ in the network keeps track of two local variables $\bbx_i, \bbd_i \in \reals^p$ which are iteratively updated at each round $t$ using the information gathered from the neighboring nodes. The vector $\bbx_i^t$ is the local decision variable of node $i$ at step $t$ whose value we expect to eventually converge to the $(1-1/e)$ fraction of the optimal solution of Problem \eqref{original_optimization_problem1}. The vector $\bbd_i^t$ is the estimate of the gradient of the global objective function that node $i$ keeps at step $t$. 
 
 To properly incorporate the received information from their neighbors, nodes should assign nonnegative weights to their neighbors. Define $w_{ij}\geq0$ to be the  weight that node $i$ assigns to node $j$. These weights indicate the effect of (variable or gradient) information nodes received from their neighbors in order to update their local (variable or gradient) information. Indeed, the weights $w_{ij}$ must fulfill some requirements (later described in Assumption \ref{ass:weights}), but they are design parameters of DCG and can be properly chosen by the nodes prior to the implementation of the algorithm.
 
The first step at each round $t$ of DCG is updating the local gradient approximation vectors $\bbd_i^t$ using local and neighboring gradient information. In particular, node $i$ computes its vector $\bbd_i^t$ according to the update rule
%%%%
\begin{align}\label{eq:gradient_update}
\bbd_i^t = (1-\alpha) \!\!\!\sum_{j\in \ccalN_i \cup\{i\}}\!\!  w_{ij}\bbd_{j}^{t-1}\ +\ \alpha \nabla F_i(\bbx_i^t),
\end{align}
%%%%
where $\alpha\in[0,1]$ is an averaging coefficient. Note that the sum  $\sum_{j\in \ccalN_i\cup\{i\}}  w_{ij}\bbd_{j}^{t-1}$ in \eqref{eq:gradient_update} is a weighted average of node $i$'s vector $\bbd_i^{t-1}$ and its neighbors $\bbd_j^{t-1}$, evaluated at step $t-1$. Hence, node $i$ computes the vector $\bbd_i^t$ by evaluating a weighted average of its current local gradient $\nabla F_i(\bbx_i^t)$ and  the local and neighboring gradient information at step $t-1$, i.e., $\sum_{j\in \ccalN_i\cup\{i\}}  w_{ij}\bbd_{j}^{t-1}$. Since the vector $\bbd_i^t$ is evaluated by aggregating gradient information from neighboring nodes, it is reasonable to expect that $\bbd_i^t$ becomes a proper approximation for the global objective function gradient $(1/n) \sum_{k=1}^n\nabla f_{k}(x)$ as time progresses. Note that to implement the update in \eqref{eq:gradient_update} nodes should exchange their local vectors $\bbd_i^t$ with their neighbors.

Using the gradient approximation vector $\bbd_i^t$, each node $i$ evaluates its local ascent direction $\bbv_i^t$ by solving the following linear program
%%%
\begin{align}\label{eq:descent_update}
\bbv_i^t = \argmax_{\bbv\in \ccalC}\ \langle \bbd_i^t, \bbv\rangle.
\end{align}
%%%
The update in \eqref{eq:descent_update} is also known as \textit{conditional gradient} update. Ideally, in a conditional gradient method, we should choose the feasible direction $\bbv\in \ccalC$ that maximizes the inner product by the full gradient vector $\frac{1}{n} \sum_{k=1}^n\nabla F_{k}(\bbx_i^t)$.  However, since in the decentralized setting the exact gradient $\frac{1}{n} \sum_{k=1}^n\nabla F_{k}(\bbx_i^t)$ is not available at the $i$-th node, we replace it by its current approximation $\bbd_i^t$ and hence we obtain the update rule \eqref{eq:descent_update}. %Indeed, this step can also be implemented in a distributed fashion since the vector $\bbd_i^t$ is locally available at node $i$.

%%%%%%%%%%%%%%%%%%%%%%%%%%%%%%%%%%%
%%%%%%%%%%%%%%%%%%%%%%%%%%%%%%%%%%%
%%%   A   L   G   O   R   I   T   H   M    %%%%%%%%%%%%%%%%
%%%%%%%%%%%%%%%%%%%%%%%%%%%%%%%%%%%
%%%%%%%%%%%%%%%%%%%%%%%%%%%%%%%%%%%
%
\begin{algorithm}[tb]
\caption{DCG at node $i$}\label{algo_DCG} 
\begin{algorithmic}[1] 
{\REQUIRE Stepsize $\alpha$ and weights $w_{ij}$ for $j\in\ccalN_{i}\cup\{i\}$
 \vspace{1mm}
   \STATE Initialize local vectors as $\bbx_i^0=\bbd_i^0=\bb0_p$ 
    \vspace{1mm}
   \STATE Initialize neighbor's vectors as $x_j^0=\bbd_j^0=\bb0_p$ if $j\in \ccalN_i$
    \vspace{1mm}
\FOR {$t=1,2,\ldots, T$}
 \vspace{1mm}
   \STATE  Compute $\displaystyle{\bbd_i^t = (1-\alpha) \!\!\!\! \sum_{j\in \ccalN_i \cup\{i\}}  \!\!\!\! w_{ij}\bbd_{j}^{t-1}+\alpha \nabla F_i(\bbx_i^t)}$;
   \vspace{1mm}
   \STATE Exchange $\bbd_i^t$ with neighboring nodes ${j\in \ccalN_i }$
   \vspace{1mm}
   \STATE  Evaluate $\bbv_i^t = \argmax_{\bbv\in \ccalC}\ \langle \bbd_i^t, \bbv\rangle$;
   \vspace{1mm}
   \STATE Update the variable $\displaystyle{\bbx_i^{t+1} =  \!\!\!\! \sum_{j\in \ccalN_i\cup\{i\}}  \!\!\!\!  w_{ij}\bbx_{j}^{t}+\frac{1}{T} \bbv_i^t}$;
    \vspace{1mm}
    \STATE Exchange $\bbx_i^{t+1}$ with neighboring nodes ${j\in \ccalN_i }$
     \vspace{1mm}
\ENDFOR}
\end{algorithmic}\end{algorithm}

After computing the local ascent directions $\bbv_i^t$, the nodes update their local variables $x_{i}^{t}$ by averaging their local and neighboring iterates and ascend in the direction $\bbv_i^t$ with stepsize $1/T$ where $T$ is the total number of iterations, i.e.,
%%%
\begin{align}\label{eq:variable_update}
\bbx_i^{t+1} =\! \sum_{j\in \ccalN_i\cup\{i\}} \!\!\! w_{ij}\bbx_{j}^{t}\ +\ \frac{1}{T} \bbv_i^t.
\end{align}
%%%
The update rule \eqref{eq:variable_update} ensures that the neighboring iterates are not far from each other via the averaging term $\sum_{j\in \ccalN_i\cup\{i\}}  w_{ij}\bbx_{j}^{t}$, while the iterates approach the optimal maximizer of the global objective function by ascending in the conditional gradient direction $\bbv_i^t$. The update in \eqref{eq:variable_update} requires a round of local communication among neighbors to exchange their local variables $\bbx_i^t$. The steps of the DCG method are summarized in Algorithm \ref{algo_DCG}.

Indeed, the weights $w_{ij}$ that nodes assign to each other cannot be arbitrary. In the following, we formalize the conditions that they should satisfy \citep{yuan2016convergence}. 
%
%
% and should satisfy some conditions. We formalize these conditions in the following \citep{yuan2016convergence}. 
%
%

%%%%%%%%%%%%%%%%%%%%%%%%%%%%%%%%%%%%%%%%%%
%%%%%%%%%%%%%%%%%%%%%%%%%%%%%%%%%%%%%%%%%%
%%%%    A  S  S  U  M  P  T  I  O  N   %%%%%%%%%%%%%%%%%%%%%%%
%%%%%%%%%%%%%%%%%%%%%%%%%%%%%%%%%%%%%%%%%%
%%%%%%%%%%%%%%%%%%%%%%%%%%%%%%%%%%%%%%%%%%
\begin{assumption}\label{ass:weights}
The weights that nodes assign to each other are nonegative, i.e., $w_{ij}\geq 0$ for all $i,j\in\ccalN$, and if node $j$ is not a neighbor of node $i$ then the corresponding weight is zero, i.e., $w_{ij}=0$ if $j\notin \ccalN_i$. Further, the weight matrix $\bbW\in\reals^{n\times n}$ with entries $w_{ij}$ satisfies
%%%%
\begin{equation}\label{eqn_conditions_on_weights}
   \bbW^{\dag}=\bbW, \quad
   \bbW\bbone_n=\bbone_n, \quad
   {\rm{null}}(\bbI-\bbW)={\rm{span}}(\bbone_n).
\end{equation}
%%%
\end{assumption}

The first condition in \eqref{eqn_conditions_on_weights} ensures that the weights are symmetric, i.e., $w_{ij}=w_{ji}$. The second condition guarantees the weights that each node assigns to itself and its neighbors sum up to 1, i.e., $\sum_{j=1}^{n} w_{ij}=1$ for all $i$. Note that the condition $\bbW\bbone_n=\bbone_n$ implies that $\bbI-\bbW$ is rank deficient. Hence, the last condition in \eqref{eqn_conditions_on_weights} ensures that the rank of $\bbI-\bbW$ is exactly $n-1$. Indeed, it is possible to optimally design the weight matrix $\bbW$ to accelerate the averaging process as discussed in \citep{boyd2004fastest}, but this is not the focus of this paper. We should emphasize that $\bbW$ is not a problem parameter, and we design it prior to runnig DCG.

Notice that the stepsize $1/T$ and the conditions in Assumption~\ref{ass:weights} on the weights $w_{ij}$ are needed to ensure that the local variables $\bbx_i^t$ are in the feasible set $\ccalC$, as stated in the following proposition.

%%%%%%%%%%%%%%%%%%%%%%%%%%%%%%%%%%%%%%%%%%
%%%%%%%%%%%%%%%%%%%%%%%%%%%%%%%%%%%%%%%%%%
%%%%     P  R  O  P  O  S  I  T  I  O  N    %%%%%%%%%%%%%%%%%%%%%
%%%%%%%%%%%%%%%%%%%%%%%%%%%%%%%%%%%%%%%%%%
%%%%%%%%%%%%%%%%%%%%%%%%%%%%%%%%%%%%%%%%%%
\begin{proposition} \label{prop_stay}
Consider the proposed DCG method outlined in Algorithm \ref{algo_DCG}. If Assumption \ref{ass:weights} holds and nodes start from $\bbx_i^0=\bb0_p\in \ccalC$, then the local iterates $\bbx_i^t$ are always in the feasible set $\ccalC$, i.e., $\bbx_i^t\in \ccalC$ for all $i\in \ccalN$ and $t=1,\dots,T$.
\end{proposition}
\begin{myproof}
Check Section \ref{app:prop_stay} in the Appendix.
\end{myproof}

Let us now explain how DCG relates to and innovates beyond the exisiting work in submodular maximization as well as decentralized convex optimization. Note that in order to solve Problem~\eqref{original_optimization_problem1} in a \emph{centralized} fashion (i.e., when every node has access to \emph{all} the local functions) we can use the continuous greedy algorithm  \citep{vondrak2008optimal}, a variant of the conditional gradient method. However, in decentralized settings, nodes have only access to their local gradients, and therefore, continuous greedy is not implementable. Similar to the decentralized convex  optimization, we can address this issue via local information aggregation. Our proposed DCG method incorporates the idea of choosing the ascent direction according to a conditional gradient update as is done in the continuous greedy algorithm (i.e., the update rule  \eqref{eq:descent_update}), while it aggregates the global objective function information through local communications with neighboring nodes (i.e., the update rule \eqref{eq:variable_update}). Unlike traditional consensus optimization methods that require exchanging nodes' local variables only  \citep{nedic2009distributed,nedic2010constrained}, DCG also requires exchanging local gradient vectors to achieve a $(1-1/e)$ fraction of the optimal solution at each node (i.e., the update rule \eqref{eq:gradient_update}). This major difference is due to the fact that in conditional gradient methods, unlike proximal gradient algorithms, the local gradients can not be used instead of the global gradient. In other words, in the update rule \eqref{eq:descent_update}, we can not use the local gradients $\nabla F_i(\bbx_i^t)$ in lieu of $\bbd_i^t$. Indeed, there are settings  for which such a replacement provides arbitrarily bad solutions.  We formally characterize the convergence of  DCG in Theorem \ref{theorem:main_theorem}.

\subsection{Extension to the Discrete Setting}

%%%%%%%%%%%%%%%%%%%%%%%%%%%%%%%%%%%%%
%%%%%%%%%%%%%%%%%%%%%%%%%%%%%%%%%%%%%
%%%%%%      S  U  B  - -  S  E  C  T  I  O  N     %%%%%%%%%%%%%
%%%%%%%%%%%%%%%%%%%%%%%%%%%%%%%%%%%%%
%%%%%%%%%%%%%%%%%%%%%%%%%%%%%%%%%%%%%
%%%%%%%%%%%%%%%%%%%%%%%%%%%%%%%%%%%%%
In this section we show how DCG can be used
 for maximizing a decentralized submodular \emph{set} function $f$, namely Problem~\eqref{original_optimization_problem2}, through its continuous relaxation. 
% 
% 
%  the multilinear extension of the function $f$. 
%  
%  
Formally, in lieu of solving Problem~\eqref{original_optimization_problem2}, we can form the following decentralized continuous optimization problem
\begin{align}\label{eq:multilinear_program}
\max_{\bbx \in \ccalC}  \frac{1}{n} \sum_{i=1}^n F_i(\bbx),
\end{align}
where $F_i$ is the multilinear extension of $f_i$ defined as 
%%%
\begin{equation}\label{eq:def_multi_linear_extension}
F_i(\bbx) = \sum_{S\subset V}f_i(S) \prod_{i\in S} x_i \prod_{j\notin S} (1-x_j) ,
\end{equation}
%%%
and the down-closed convex set $\ccalC= \text{conv}\{1_{I} : I\in \ccalI \}$ is the matroid polytope. Note that the discrete and continuous optimization formulations lead to the same optimal value \citep{calinescu2011maximizing}. 

Based on the expression in \eqref{eq:def_multi_linear_extension}, computing the full gradient $\nabla F_i$ at each node $i$ will require an exponential computation in terms of $|V|$, since the number of summands in \eqref{eq:def_multi_linear_extension} is $2^{|V|}$. As a result, in the discrete setting, we will slightly modify the DCG algorithm and work with \emph{unbiased estimates} of the gradient that can be computed in  time $O(|V|)$ (see Appendix~\ref{unbiased} for one such estimator).  More precisely, in the discrete setting, each node $i \in \nodes$ updates three local variables $\bbx_i^t, \bbd_i^t, \bbg_i^t \in \reals^{|V|}$. The variables $\bbx_i^t, \bbd_i^t$ play the same role as in DCG and are  updated using the messages received from the neighboring nodes. The variable $\bbg_i^t$ at node $i$ is defined to approximate the local gradient $\nabla F_i (\bbx_i^t)$. Consider the vector $\nabla \tilde{F}_i(\bbx_i^t)$ as an unbiased estimator of the local gradient $\nabla F_i (\bbx_i^t)$ at time $t$, and define the vector $\bbg_i^t$ as the outcome of the recursion 
%%%
\begin{equation}\label{gradient_approx_update}
\bbg_i^t = (1-\phi) \bbg_i^{t-1}+\phi \nabla \tilde{F}_i(\bbx_i^t),
\end{equation}
%%%
where $\phi\in[0,1]$ is the averaging parameter. We initialize all vectors as $\bbg_i^{0}=\bb0\in \reals^{|V|}$. It was shown recently \citep{mokhtari2017conditional} that  the averaging technique in \eqref{gradient_approx_update} reduces the noise of the gradient approximations. Therefore,  the sequence of $\bbg_i^t$ approaches the true local gradient $ \nabla {F}_i(\bbx_i^t)$ as time progresses. 

%%%%%%%%%%%%%%%%%%%%%%%%%%%%%%%%%%%
%%%%%%%%%%%%%%%%%%%%%%%%%%%%%%%%%%%
%%%   A   L   G   O   R   I   T   H   M    %%%%%%%%%%%%%%%%
%%%%%%%%%%%%%%%%%%%%%%%%%%%%%%%%%%%
%%%%%%%%%%%%%%%%%%%%%%%%%%%%%%%%%%%
%
\begin{algorithm}[tb]
\caption{Discrete DCG at node $i$}\label{algo_DDCG} 
\begin{algorithmic}[1] 
{\REQUIRE $\alpha,\phi\in[0,1]$ and weights $w_{ij}$ for $j\in\ccalN_{i}\cup\{i\}$;
 \vspace{1mm}
   \STATE Initialize local vectors as $\bbx_i^0=\bbd_i^0=\bbg_i^0=\bb0$ ;
    \vspace{1mm}
   \STATE Initialize neighbor's vectors as $\bbx_j^0=\bbd_j^0=\bb0$ if $j\in \ccalN_i$;
    \vspace{1mm}
\FOR {$t=1,2,\ldots, T$}
 \vspace{1mm}
   \STATE  Compute $\displaystyle{\bbg_i^t = (1-\phi) \bbg_i^{t-1}+\phi \nabla \tilde{F}_i(\bbx_i^t)}$;
    \vspace{1mm}
   \STATE  Compute $\displaystyle{\bbd_i^t = (1-\alpha) \!\!\!\! \sum_{j\in \ccalN_i \cup\{i\}}  \!\!\!\! w_{ij}\bbd_{j}^{t-1}+\alpha \bbg_i^t}$;
   \vspace{1mm}
   \STATE Exchange $\bbd_i^t$ with neighboring nodes ${j\in \ccalN_i }$;
   \vspace{1mm}
   \STATE  Evaluate $\bbv_i^t = \argmax_{\bbv\in \ccalC}\ \langle \bbd_i^t, \bbv\rangle$;
   \vspace{1mm}
   \STATE Update the variable $\displaystyle{\bbx_i^{t+1} =  \!\!\!\! \sum_{j\in \ccalN_i\cup\{i\}}  \!\!\!\!  w_{ij}\bbx_{j}^{t}+\frac{1}{T} \bbv_i^t}$;
    \vspace{1mm}
    \STATE Exchange $\bbx_i^{t+1}$ with neighboring nodes ${j\in \ccalN_i }$;
     \vspace{1mm}
\ENDFOR}
\STATE Apply proper rounding to obtain a solution for \eqref{original_optimization_problem2};
\end{algorithmic}\end{algorithm}

The steps of the \alg for the discrete setting is summarized in Algorithm~\ref{algo_DDCG}. Note that the major difference between the Discrete DCG  method (Algorithm~\ref{algo_DDCG})  and the continuous DCG method (Algorithm~\ref{algo_DCG}) is in Step 5 in which the exact local gradient $\nabla F_i (\bbx_i^t)$ is replaced by the stochastic approximation $\bbg_i^t$ which only requires access to the computationally cheap unbiased gradient estimator $\nabla \tilde{F}_i(\bbx_i^t)$. 
The communication complexity of both the discrete and continuous versions of DCG are the same at each round. However, since we are using unbiased estimations of the local gradients $\nabla F_i(\bbx_i)$, the  Discrete DCG takes more  rounds to converge to a near-optimal solution compared to~continuous DCG. We characterize the convergence of Discrete DCG in Theorem \ref{theorem:main_theorem_discrete}. {Further, the implementation of Discrete DCG requires rounding the continuous solution to obtain a discrete solution for the original problem without any loss in terms of the objective function value. The provably lossless rounding schemes include the pipage rounding \citep{calinescu2011maximizing} and contention resolution~\citep{DBLP:journals/siamcomp/ChekuriVZ14}.
}

\section{Convergence Analysis} \label{sec:convergence}

In this section, we study the convergence properties of DCG in both continuous and discrete settings. In this regard, we assume that the following conditions hold.

%%%%%%%%%%%%%%%%%%%%%%%%%%%%%%%%%%%%%%%%%%
%%%%%%%%%%%%%%%%%%%%%%%%%%%%%%%%%%%%%%%%%%
%%%%    A  S  S  U  M  P  T  I  O  N   %%%%%%%%%%%%%%%%%%%%%%%
%%%%%%%%%%%%%%%%%%%%%%%%%%%%%%%%%%%%%%%%%%
%%%%%%%%%%%%%%%%%%%%%%%%%%%%%%%%%%%%%%%%%%
\begin{assumption}\label{ass:bounded_set}
{Euclidean distance of the elements in the set $\ccalC$ are uniformly bounded, i.e., for all $\bbx,\bby \in \ccalC$ we have}
%%%
\begin{equation}
\|\bbx-\bby\|\leq D.
\end{equation}
%%%
\end{assumption}
%%%%%%%%%%%%%%%%%%%%%%%%%%%%%%%%%%%%%%%%%%
%%%%%%%%%%%%%%%%%%%%%%%%%%%%%%%%%%%%%%%%%%
%%%%    A  S  S  U  M  P  T  I  O  N   %%%%%%%%%%%%%%%%%%%%%%%
%%%%%%%%%%%%%%%%%%%%%%%%%%%%%%%%%%%%%%%%%%
%%%%%%%%%%%%%%%%%%%%%%%%%%%%%%%%%%%%%%%%%%
\begin{assumption}\label{ass:smoothness}
The local objective functions $F_i(\bbx)$ are monotone and DR-submodular. Further, their gradients are $L$-Lipschitz continuous over the set $\ccalX$, i.e., for all $\bbx,\bby \in \ccalX$ 
%%%
\begin{equation}
\| \nabla F_i(\bbx) -  \nabla F_i(\bby) \| \leq L \| \bbx-\bby \|.
\end{equation}
%%%
\end{assumption}

%%%%%%%%%%%%%%%%%%%%%%%%%%%%%%%%%%%%%%%%%%
%%%%%%%%%%%%%%%%%%%%%%%%%%%%%%%%%%%%%%%%%%
%%%%    A  S  S  U  M  P  T  I  O  N   %%%%%%%%%%%%%%%%%%%%%%%
%%%%%%%%%%%%%%%%%%%%%%%%%%%%%%%%%%%%%%%%%%
%%%%%%%%%%%%%%%%%%%%%%%%%%%%%%%%%%%%%%%%%%
\begin{assumption}\label{ass:smoothness2}
The norm of gradients  $\|\nabla F_i(\bbx)\|$ are bounded over the convex set $\ccalC$, i.e., for all $\bbx \in \ccalC$, $i\in\ccalN$,
%%%
\begin{equation}
\| \nabla F_i(\bbx) \| \leq G.
\end{equation}
%%%
\end{assumption}

The condition in Assumption \ref{ass:bounded_set} guarantees that the diameter of the convex set $\ccalC$ is bounded. Assumption \ref{ass:smoothness} is needed to ensure that the local objective functions $F_i$ are smooth. Finally, the condition in Assumption \ref{ass:smoothness2} enforces the gradients norm to be bounded over the convex set $\ccalC$. All these assumptions are customary and necessary in the analysis of decentralized algorithms. For more details, please check Section VII-B in \citet{jakovetic2014fast}.

%%%%%%%%%%%%%%%%%%%%%%%%%%%%%%%%%%%%%%%%%%
%%%%%%%%%%%%%%%%%%%%%%%%%%%%%%%%%%%%%%%%%%
%%%%    M  A  I  N    M  A  T  T  E  R    %%%%%%%%%%%%%%%%%%%%%%
%%%%%%%%%%%%%%%%%%%%%%%%%%%%%%%%%%%%%%%%%%
%%%%%%%%%%%%%%%%%%%%%%%%%%%%%%%%%%%%%%%%%

We proceed to derive a constant factor approximation for DCG. Our main result is stated in Theorem~\ref{theorem:main_theorem}. However, to better illustrate the main result, we first need to provide several definitions and technical lemmas.  Let us begin by defining the average variables $\bar{\bbx}^t$ as
\begin{equation} \label{average_x}
\bar{\bbx}^{t}= \frac{1}{n}\sum_{i=1}^n \bbx_i^t.
\end{equation}

In the following lemma, we establish an upper bound on the variation in the sequence of average variables $\{\bar{\bbx}^t\}$.

%%%%%%%%%%%%%%%%%%%%%%%%%%%%%%%%%%%%%%%%%%
%%%%%%%%%%%%%%%%%%%%%%%%%%%%%%%%%%%%%%%%%%
%%%%    L   E   M   M   A    %%%%%%%%%%%%%%%%%%%%%%%%%%%%
%%%%%%%%%%%%%%%%%%%%%%%%%%%%%%%%%%%%%%%%%%
%%%%%%%%%%%%%%%%%%%%%%%%%%%%%%%%%%%%%%%%%%
\begin{lemma}\label{lemma:ar_in_avg_bound}
Consider the proposed DCG algorithm defined in Algorithm \ref{algo_DCG}. Further, recall the definition of $\bar{\bbx}^{t}$ in \eqref{average_x}. If Assumptions \ref{ass:weights} and \ref{ass:bounded_set} hold, then the difference between two consecutive average vectors is upper bounded by
\begin{align}\label{eq:var_in_avg_bound}
\|\bar{\bbx}^{t+1} - \bar{\bbx}^{t}\| \leq \frac{D}{T} .
\end{align}
\end{lemma}
%%%%%%%%%%%%%%%%%%%%%%%%%%%%%%%%%%%%%%%%%%
%%%%%%%%%%%%%%%%%%%%%%%%%%%%%%%%%%%%%%%%%%
%%%%      P   R   O   O   F     %%%%%%%%%%%%%%%%%%%%%%%%%%%
%%%%%%%%%%%%%%%%%%%%%%%%%%%%%%%%%%%%%%%%%%
%%%%%%%%%%%%%%%%%%%%%%%%%%%%%%%%%%%%%%%%%%
\begin{myproof}
Check Section \ref{app:lemma:ar_in_avg_bound} in the Appendix.
\end{myproof}

%%%%%%%%%%%%%%%%%%%%%%%%%%%%%%%%%%%%%%%%%%
%%%%%%%%%%%%%%%%%%%%%%%%%%%%%%%%%%%%%%%%%%
%%%%    M  A  I  N    M  A  T  T  E  R    %%%%%%%%%%%%%%%%%%%%%%
%%%%%%%%%%%%%%%%%%%%%%%%%%%%%%%%%%%%%%%%%%
%%%%%%%%%%%%%%%%%%%%%%%%%%%%%%%%%%%%%%%%%
Recall that at every node $i$, the messages are mixed using the coefficients $w_{ij}$, i.e., the $i$-th row of the matrix $\bbW$.
It is thus not hard to see that the spectral properties of $\bbW$ (e.g. the spectral gap) play an important role in the 
 the speed of achieving consensus in decentralized methods.
%%%%%%%%%%%%%%%%%%%%%%%%%%%%%%%%%%%%%%%%%%
%%%%%%%%%%%%%%%%%%%%%%%%%%%%%%%%%%%%%%%%%%
%%%%     D  E  F  I  N  I  T  I  O  N     %%%%%%%%%%%%%%%%%%%%%%%
%%%%%%%%%%%%%%%%%%%%%%%%%%%%%%%%%%%%%%%%%%
%%%%%%%%%%%%%%%%%%%%%%%%%%%%%%%%%%%%%%%%%%
\begin{definition} 
Consider the eigenvalues of $\bbW$ which can be sorted in a nonincreasing order as $1 = \lambda_{1}(\bbW) \geq \lambda_{2}(\bbW) \dots  \geq  \lambda_{n}(\bbW) > -1$. Define $\beta$ as the second largest magnitude of the eigenvalues of
$\bbW$, i.e., 
%%%
\begin{align}\label{eq:def_beta}
\beta:= \max\{ | \lambda_{2}(\bbW) |,  | \lambda_{n}(\bbW) | \}.
\end{align}
\end{definition}
 As we will see, a mixing matrix $\bbW$ with smaller $\beta$ has a larger spectral gap $1-\beta$ which yields faster convergence \citep{boyd2004fastest,duchi2012dual}. In the following lemma, we derive an upper bound on the sum of the distances between the local iterates $\bbx_{i}^t$ and their average $\bar{\bbx}^{t}$, where the bound is a function of the graph spectral gap $1-\beta$, size of the network $n$, and the total number of iterations $T$.

%%%%%%%%%%%%%%%%%%%%%%%%%%%%%%%%%%%%%%%%%%
%%%%%%%%%%%%%%%%%%%%%%%%%%%%%%%%%%%%%%%%%%
%%%%    L   E   M   M   A    %%%%%%%%%%%%%%%%%%%%%%%%%%%%
%%%%%%%%%%%%%%%%%%%%%%%%%%%%%%%%%%%%%%%%%%
%%%%%%%%%%%%%%%%%%%%%%%%%%%%%%%%%%%%%%%%%%
\begin{lemma}\label{lemma:eq:bound_on_dif_from_avg}
Consider the proposed DCG algorithm defined in Algorithm \ref{algo_DCG}. Further, recall the definition of $\bar{\bbx}^{t}$ in \eqref{average_x}. If Assumptions \ref{ass:weights} and \ref{ass:bounded_set} hold, then for all  $t\leq T$ we have %the sum of squared Euclidean norm of the difference between the local iterates $\bbx_i^t$ and their average $\bar{\bbx}^{t}$ is upper bounded by
%%%
\begin{align}\label{eq:bound_on_dif_from_avg}
\left( \sum_{i=1}^n \left\|\bbx_i^t-\bar{\bbx}^t\right\|^2\right)^{1/2} \leq  \frac{\sqrt{n}D}{T(1-\beta)}.
\end{align}
%for any $t\leq T$.
%%%
\end{lemma}
%%%%%%%%%%%%%%%%%%%%%%%%%%%%%%%%%%%%%%%%%%
%%%%%%%%%%%%%%%%%%%%%%%%%%%%%%%%%%%%%%%%%%
%%%%      P   R   O   O   F     %%%%%%%%%%%%%%%%%%%%%%%%%%%
%%%%%%%%%%%%%%%%%%%%%%%%%%%%%%%%%%%%%%%%%%
%%%%%%%%%%%%%%%%%%%%%%%%%%%%%%%%%%%%%%%%%%
\begin{myproof}
Check Section \ref{app:lemma:eq:bound_on_dif_from_avg} in the Appendix. 
\end{myproof}

%%%%%%%%%%%%%%%%%%%%%%%%%%%%%%%%%%%%%%%%%%
%%%%%%%%%%%%%%%%%%%%%%%%%%%%%%%%%%%%%%%%%%
%%%%    M  A  I  N    M  A  T  T  E  R    %%%%%%%%%%%%%%%%%%%%%%
%%%%%%%%%%%%%%%%%%%%%%%%%%%%%%%%%%%%%%%%%%
%%%%%%%%%%%%%%%%%%%%%%%%%%%%%%%%%%%%%%%%%
Let us now define $\bar{\bbd}^t $ as the average of local gradient approximations $\bbd_i^t$ at step $t$, i.e., 
\begin{equation}
\bar{\bbd}^t = \frac{1}{n} \sum_{i=1}^n \bbd_i^t.
\end{equation}
We will show in the following that the vectors $\bbd_i^t$ also become uniformly close to $\bar{\bbd}^t$.

%%%%%%%%%%%%%%%%%%%%%%%%%%%%%%%%%%%%%%%%%%
%%%%%%%%%%%%%%%%%%%%%%%%%%%%%%%%%%%%%%%%%%
%%%%    L   E   M   M   A    %%%%%%%%%%%%%%%%%%%%%%%%%%%%
%%%%%%%%%%%%%%%%%%%%%%%%%%%%%%%%%%%%%%%%%%
%%%%%%%%%%%%%%%%%%%%%%%%%%%%%%%%%%%%%%%%%%
\begin{lemma}\label{lemma:bound_on_gradient_consensus_error}
Consider the proposed DCG algorithm defined in Algorithm \ref{algo_DCG}. If Assumptions \ref{ass:weights} and \ref{ass:smoothness} hold, then %Euclidean norm of the difference of between the  approximate gradient $\bbd_i^t $ and the average of approximate gradients  $\bar{\bbd}^t$ at step $t$ is upper bounded by 
%%%
\begin{align}\label{eq:bound_on_gradient_consensus_error}
\left(\sum_{i=1}^n\|\bbd_i^t-\bar{\bbd}^t\|^2\right)^{1/2} \leq \frac{\alpha \sqrt{n} G}{1-\beta(1-\alpha)}.
\end{align}
%%%
\end{lemma}
%%%%%%%%%%%%%%%%%%%%%%%%%%%%%%%%%%%%%%%%%%
%%%%%%%%%%%%%%%%%%%%%%%%%%%%%%%%%%%%%%%%%%
%%%%      P   R   O   O   F     %%%%%%%%%%%%%%%%%%%%%%%%%%%
%%%%%%%%%%%%%%%%%%%%%%%%%%%%%%%%%%%%%%%%%%
%%%%%%%%%%%%%%%%%%%%%%%%%%%%%%%%%%%%%%%%%%
\begin{myproof}
Check Section \ref{proof:lemma:bound_on_gradient_consensus_error} in the Appendix. 
\end{myproof}

%%%%%%%%%%%%%%%%%%%%%%%%%%%%%%%%%%%%%%%%%%
%%%%%%%%%%%%%%%%%%%%%%%%%%%%%%%%%%%%%%%%%%
%%%%    M  A  I  N    M  A  T  T  E  R    %%%%%%%%%%%%%%%%%%%%%%
%%%%%%%%%%%%%%%%%%%%%%%%%%%%%%%%%%%%%%%%%%
%%%%%%%%%%%%%%%%%%%%%%%%%%%%%%%%%%%%%%%%%
Lemma \ref{lemma:bound_on_gradient_consensus_error} guarantees that the individual local gradient approximation vectors $\bbd_i^t$ are close to the average vector $\bar{\bbd}^t$ if the parameter $\alpha$ is small. To show that the gradient vectors $\bbd_i^t$, generated by DCG, approximate the gradient of the global objective function, we further need to show that the average vector  $\bar{\bbd}^t$ approaches the global objective function gradient $\nabla F$. We  prove this claim in the following lemma. 

%%%%%%%%%%%%%%%%%%%%%%%%%%%%%%%%%%%%%%%%%%
%%%%%%%%%%%%%%%%%%%%%%%%%%%%%%%%%%%%%%%%%%
%%%%    L   E   M   M   A    %%%%%%%%%%%%%%%%%%%%%%%%%%%%
%%%%%%%%%%%%%%%%%%%%%%%%%%%%%%%%%%%%%%%%%%
%%%%%%%%%%%%%%%%%%%%%%%%%%%%%%%%%%%%%%%%%%
\begin{lemma}\label{lemma:bound_on_gradient_error}
Consider the proposed DCG algorithm defined in Algorithm \ref{algo_DCG}. If Assumptions \ref{ass:weights}-\ref{ass:smoothness2} hold, then
%%%
\begin{align}\label{eq:bound_on_gradient_error}
\left\|\bar{\bbd}^t   - \frac{1}{n}\sum_{i=1}^n\nabla F_i(\bar{\bbx}^t) \right\|
\leq (1-\alpha)^t G
+\left(\frac{(1-\alpha)LD}{\alpha T} 
+ \frac{ LD}{T(1-\beta)} \right).
\end{align}
%%%
\end{lemma}
%%%%%%%%%%%%%%%%%%%%%%%%%%%%%%%%%%%%%%%%%%
%%%%%%%%%%%%%%%%%%%%%%%%%%%%%%%%%%%%%%%%%%
%%%%      P   R   O   O   F     %%%%%%%%%%%%%%%%%%%%%%%%%%%
%%%%%%%%%%%%%%%%%%%%%%%%%%%%%%%%%%%%%%%%%%
%%%%%%%%%%%%%%%%%%%%%%%%%%%%%%%%%%%%%%%%%
\begin{myproof}
Check Section \ref{proof:lemma:bound_on_gradient_error} in the Appendix.
\end{myproof}

%%%%%%%%%%%%%%%%%%%%%%%%%%%%%%%%%%%%%%%%%%
%%%%%%%%%%%%%%%%%%%%%%%%%%%%%%%%%%%%%%%%%%
%%%%    M  A  I  N    M  A  T  T  E  R    %%%%%%%%%%%%%%%%%%%%%%
%%%%%%%%%%%%%%%%%%%%%%%%%%%%%%%%%%%%%%%%%%
%%%%%%%%%%%%%%%%%%%%%%%%%%%%%%%%%%%%%%%%%
By combining Lemmas \ref{lemma:bound_on_gradient_consensus_error} and \ref{lemma:bound_on_gradient_error} and setting $\alpha=1/\sqrt{T}$ we can conclude that the local gradient approximation vector $\bbd_i^t$ of each node $i$ is within $\mathcal{O}(1/\sqrt{T})$ distance of the global objective gradient $\nabla F(\bar{\bbx}^t)$ evaluated at $\bar{\bbx}^t$. We use this observation in the following theorem to show that the sequence of iterates generated by DCG achieves the tight $(1-1/e)$ approximation ratio of the optimum global solution.

%%%%%%%%%%%%%%%%%%%%%%%%%%%%%%%%%%%%%%%%%%
%%%%%%%%%%%%%%%%%%%%%%%%%%%%%%%%%%%%%%%%%%
%%%%    T  H  E  O  R  E  M    %%%%%%%%%%%%%%%%%%%%%%%%%%
%%%%%%%%%%%%%%%%%%%%%%%%%%%%%%%%%%%%%%%%%%
%%%%%%%%%%%%%%%%%%%%%%%%%%%%%%%%%%%%%%%%%%
\begin{theorem}\label{theorem:main_theorem}
Consider the proposed DCG method outlined in Algorithm~\ref{algo_DCG}. Further, consider $\bbx^*$ as the global maximizer of Problem \eqref{original_optimization_problem1}. If Assumptions \ref{ass:weights}-\ref{ass:smoothness2} hold and we set $\alpha= 1/\sqrt{T}$, %the variables average $\bar{\bbx}^T$ at step $T$ satisfies 
%%%%
%\begin{align}\label{eq:main_result_1}
%F(\bar{\bbx}^T)
%&\geq (1-e^{-1} ) F(\bbx^*)-  \frac{ LD^2+GD}{T^{1/2}}-\frac{GD}{T^{1/2}(1-\beta)}
%\nonumber\\
%&
%\quad - \frac{ LD^2}{2T} - \frac{ LD^2}{T(1-\beta)} ,
%\end{align}
%%%%
%where $\bbx^*$ is the global maximizer of Problem \eqref{original_optimization_problem1}. \red{Should we keep \eqref{eq:main_result_1}?} Moreover, 
for all nodes $j\in \ccalN$, the local variable $\bbx_j^T$ obtained after $T$ iterations satisfies 
\begin{align}\label{local_node_bound}
F(\bbx_j^T)&\geq (1-e^{-1} )F(\bbx^*) - \frac{ LD^2+GD}{T^{1/2}}-\frac{GD}{T^{1/2}(1-\beta)}-  \frac{ LD^2}{2T} - \frac{ GD+LD^2}{T(1-\beta)}.
\end{align}
\end{theorem}
%%%%%%%%%%%%%%%%%%%%%%%%%%%%%%%%%%%%%%%%%%
%%%%%%%%%%%%%%%%%%%%%%%%%%%%%%%%%%%%%%%%%%
%%%%      P   R   O   O   F     %%%%%%%%%%%%%%%%%%%%%%%%%%%
%%%%%%%%%%%%%%%%%%%%%%%%%%%%%%%%%%%%%%%%%%
%%%%%%%%%%%%%%%%%%%%%%%%%%%%%%%%%%%%%%%%%%
\begin{myproof}
Check Section \ref{proof:theorem:main_theorem} in the Appendix.
\end{myproof}

Theorem \ref{theorem:main_theorem} shows that the sequence of the local variables ${\bbx}_j^t$, generated by DCG, is able to achieve the optimal approximation ratio $(1-1/e)$, while the error term vanishes at a
sublinear rate of $\mathcal{O}(1/T^{1/2})$, i.e., 
%%%
\begin{align}
F(\bbx_j^T)
\geq (1-1/e ) F(\bbx^*)-  \mathcal{O}\left(\frac{1}{(1-\beta)T^{1/2}}\right),
\end{align}
%%%
which implies that the iterate of \textit{each node} reaches an objective value larger than $(1-1/e-\epsilon)OPT$ after $\mathcal{O}(1/\eps^2)$ rounds of communication. It is worth mentioning that the result in Theorem \ref{theorem:main_theorem} is consistent with classical results in decentralized optimization that the error term vanishes faster for the graphs with larger  spectral gap $1-\beta$. We proceed to study the convergence properties of Discrete DCG in Algorithm \ref{algo_DDCG}. To do so, we first assume that the variance of the stochastic gradients $  \nabla \tilde{F}_i(\bbx)$ used in Discrete DCG is bounded. We justify this assumption in Remark~\ref{remdec}.

%%%%%%%%%%%%%%%%%%%%%%%%%%%%%%%%%%%%%%%%%%
%%%%%%%%%%%%%%%%%%%%%%%%%%%%%%%%%%%%%%%%%%
%%%%    A  S  S  U  M  P  T  I  O  N   %%%%%%%%%%%%%%%%%%%%%%%
%%%%%%%%%%%%%%%%%%%%%%%%%%%%%%%%%%%%%%%%%%
%%%%%%%%%%%%%%%%%%%%%%%%%%%%%%%%%%%%%%%%%%
\begin{assumption}\label{ass:bounded_variance}
The variance of the unbiased estimators $\nabla \tilde{F}(\bbx)$ is bounded above by $\sigma^2$ over the convex set $\ccalC$, i.e., for any $i \in \nodes$ and any vector $\bbx\in\ccalC$ we can write 
\begin{equation}
\E{\|  \nabla \tilde{F}_i(\bbx) - \nabla F_i(\bbx)  \|^2} \leq \sigma^2,
\end{equation}
where the expectation is with respect to the randomness of the unbiased estimator.
\end{assumption}

In the following theorem, we show that Discrete DCG achieves a $(1-1/e)$ approximation ration for Problem~\eqref{original_optimization_problem2}.

%%%%%%%%%%%%%%%%%%%%%%%%%%%%%%%%%%%%%%%%%%
%%%%    T  H  E  O  R  E  M    %%%%%%%%%%%%%%%%%%%%%%%%%%
%%%%%%%%%%%%%%%%%%%%%%%%%%%%%%%%%%%%%%%%%%
%%%%%%%%%%%%%%%%%%%%%%%%%%%%%%%%%%%%%%%%%%
\begin{theorem}\label{theorem:main_theorem_discrete}
Consider our proposed Discrete DCG algorithm outlined in Algorithm \ref{algo_DDCG}. Recall the definition of the multilinear extension function $F_i$ in \eqref{eq:def_multi_linear_extension}. If Assumptions \ref{ass:weights}-\ref{ass:bounded_variance} hold and we set $\alpha= T^{-1/2}$ and $\phi=T^{-2/3}$, then for all nodes $j\in \ccalN$ the local variables $\bbx_j^T$ obtained after running Discrete DCG for $T$ iterations satisfy 
%%%%
\begin{align}\label{eq:main_result_1}
\E{F(\bbx_j^T)}
&\geq (1-e^{-1} ) F(\bbx^*)
-\frac{ GD+LD^2}{T(1-\beta)}   -\frac{LD^2}{2T} -\frac{\sqrt{6} LD^2}{T^{2/3}}-\frac{\sqrt{12} LD^2}{(1-\beta)T^{2/3}}
     \nonumber\\
&  \quad  
      - \frac{ D(\sigma^2+G^2)^{1/2}}{T^{1/2}(1-\beta)}
   - \frac{DG+LD^2}{T^{1/2}}
   -\frac{\sqrt{2}\sigma+\sqrt{12} LD^2+4DG}{T^{1/3}}
      -\frac{ \sqrt{24}LD^2}{(1-\beta)T^{1/3}},
\end{align}
%%%
where $\bbx^*$ is the global maximizer of Problem \eqref{eq:multilinear_program}.
%Moreover, for all nodes $j=\ccalN$, the local variable $\bbx_j^T$ obtained after $T$ iterations satisfies 
\end{theorem}
%%%%%%%%%%%%%%%%%%%%%%%%%%%%%%%%%%%%%%%%%%
%%%%%%%%%%%%%%%%%%%%%%%%%%%%%%%%%%%%%%%%%%
%%%%      P   R   O   O   F     %%%%%%%%%%%%%%%%%%%%%%%%%%%
%%%%%%%%%%%%%%%%%%%%%%%%%%%%%%%%%%%%%%%%%%
%%%%%%%%%%%%%%%%%%%%%%%%%%%%%%%%%%%%%%%%%%
\begin{myproof}
Check Section \ref{proof:theorem:main_theorem_discrete} in the Appendix.
\end{myproof}

Theorem  \ref{theorem:main_theorem_discrete} states that the sequence of the local variables ${\bbx}_j^t$, generated by Discrete DCG, is able to achieve the optimal approximation ratio $(1-1/e)$ in expectation, while the error term vanishes at a sublinear rate of $\mathcal{O}(1/T^{1/3})$, i.e., 
%%%
\begin{equation} \label{eq:main_result_1000}
\E{F(\bbx_j^T)}\geq (1-e^{-1} )F(\bbx^*) -  \mathcal{O}\left(\frac{1}{(1\!-\!\beta)T^{1/3}}\right). 
\end{equation}
%%%%
Hence, the iterate of \textit{each node} reaches an objective value larger than $(1-1/e-\epsilon)OPT$ after $\mathcal{O}(1/\eps^3)$ rounds of communication.

%%%%%%%%%%%%%%%%%%%%%%%%%%%%%%%%%%%%%%%%%%
%%%%%%%%%%%%%%%%%%%%%%%%%%%%%%%%%%%%%%%%%%
%%%%      R. E  M. A. R. K     %%%%%%%%%%%%%%%%%%%%%%%%%%%
%%%%%%%%%%%%%%%%%%%%%%%%%%%%%%%%%%%%%%%%%%
%%%%%%%%%%%%%%%%%%%%%%%%%%%%%%%%%%%%%%%%%%
\begin{remark} \label{remdec}
For any submodular set function $h: 2^V \to \mathbb{R}$ with associated multilinear extension $H$, it can be shown that its Lipschitz constant $L$ and the gradient norm $G$ are both bounded above by 
$m_f \sqrt{|V|}$, where $m_f$ is the maximum marginal value of $f$, i.e., $m_f = \max_{i \in V} f(\{i\})$ (see, \citet{hassani2017gradient}). Similarly, it can be shown that for the unbiased estimator in Appendix~\ref{unbiased} we have $\sigma \leq m_f \sqrt{|V|}$.
\end{remark}

%!TEX root = example_paper.tex
%%%%%%%%%%%%%%%%%%%%%%%%%%%%%%%%%%%%%
%%%%%%%%%%%%%%%%%%%%%%%%%%%%%%%%%%%%%
%%%%%%      S  E  C T  I  O  N     %%%%%%%%%%%%%%%%%%%
%%%%%%%%%%%%%%%%%%%%%%%%%%%%%%%%%%%%%
%%%%%%%%%%%%%%%%%%%%%%%%%%%%%%%%%%%%%
%%%%%%%%%%%%%%%%%%%%%%%%%%%%%%%%%%%%%
\section{Numerical Experiments}
We will consider a discrete setting for our experiments and use Algorithm~\ref{algo_DDCG} to find a decentralized solution.
The main objective is to demonstrate how consensus is reached and~how the global objective increases depending on the topology of the network and the parameters of the algorithm.
% (e.g. the number of iterations).  

For our experiments, we have used the MovieLens data set. It consists of 1 million ratings (from 1 to 5) by $M=6000$ users for $p=4000$ movies.  We consider a network of $n = 100$ nodes. The data has been distributed equally between the nodes of the network, i.e., the set of users has been partitioned into $100$ equally-sized sets and each node in the network has access to only one chunk (partition) of the data. The global task is to find a set of $k$ movies that are most satisfactory to \emph{all} the users (the precise formulation will appear shortly). However, as each of the nodes in the network has access to the data of a small portion of the users, the nodes have to cooperate (exchange information) in order to fulfil the global task.

\begin{figure}[t] 
\begin{center} 
 \centerline{\includegraphics[width=.65\columnwidth]{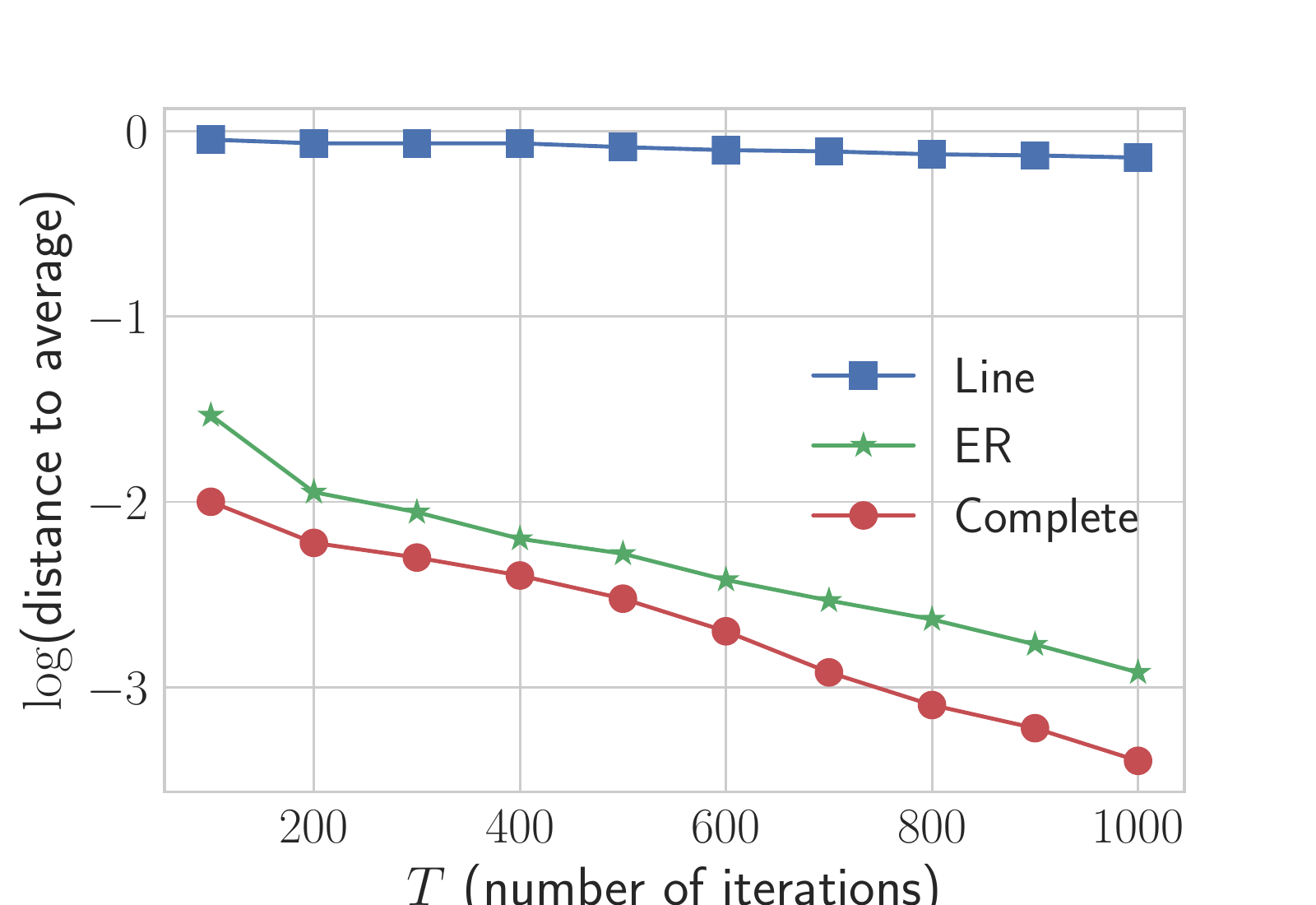}}
\end{center}
%\vspace{-6mm}
\caption{The logarithm of the distance-to-average at final round $T$ is plotted as a function of $T$. Note that when the underlying graph is complete or Erdos-Renyi (ER) with a good average degree, then consensus will be achieved even for small number of iterations $T$. However, for poor connected graphs such as the line graph, reaching consensus requires a  large number of iterations.}
\label{fig-dist}
\end{figure}

We consider a well motivated objective function for the experiments.   Let $r_{\ell,j}$ denote the rating of user $\ell$ for movie $j$ (if such a rating does not exist in the data we assign $r_{\ell,j}$ to 0). We associate to each user $\ell$ a ``facility location" objective function $g_\ell (S) = \max_{j\in S} r_{\ell,j}$, where $S$ is any subset of the movies (i.e. the ground set $V$ is the set of the movies). Such a function shows how much user $\ell$ will be ``satisfied" by a subset $S$ of the movies. Recall that each node $i$ in the network has access to the data of a (small) subset of users which we denote by $\mathcal{U}_i$. The objective function associated with node $i$ is given by $f_i(S) = \sum_{\ell \in \mathcal{U}_i} g_\ell(S)$. With such a choice of the local functions, our global task is hence to solve problem~\eqref{original_optimization_problem2} when the matroid $\mathcal{I}$ is the $k$-uniform matroid (a.k.a. the $k$-cardinality constraint).    

We consider three different choices for the underlying communication graph between the $100$ nodes: A line graph (which looks like a simple path from node 1 to node 100), an Erdos-Renyi random graph (with average degree $5$), and a complete graph.    The matrix $\bbW$ is chosen as follows (based on each of the three graphs). If $(i,j)$ is and edge of the graph, we let $w_{i,j} = 1/(1+\max(d_i, d_j))$. If $(i,j)$ is not an edge and $i,j$ are distinct integers, we have $w_{i,j} = 0$. Finally we let $w_{i,i} = 1 -\sum_{j \in \mathcal{N}} w_{i,j}$. It is not hard to show that the above choice for $\bbW$ satisfies Assumption~\ref{ass:weights}. 
 
Figure~\ref{fig-dist} shows how consensus is reached w.r.t each of the three underlying networks. To measure consensus, we plot the (logarithm of) distance-to-average value $\frac{1}{n} \sum_{i=1}^n || \bbx_i^T - \bar{\bbx}^T ||$ as a function of the total number of iterations $T$ averaged over many trials (see \eqref{average_x} for the definition of $\bar{\bbx}^T$).  It is easy to see that the distance to average is small if and only if all the local decisions $\bbx_i^T$ are close to the average decision $\bar{\bbx}^T$. As expected, it takes much less time to reach consensus when the underlying graph is fully connected (i.e. complete graph). For the line graph, the convergence is very slow as this graph has the least degree of connectivity. 

\begin{figure}[t] 
\begin{center} 
 \centerline{\includegraphics[width=.65\columnwidth]{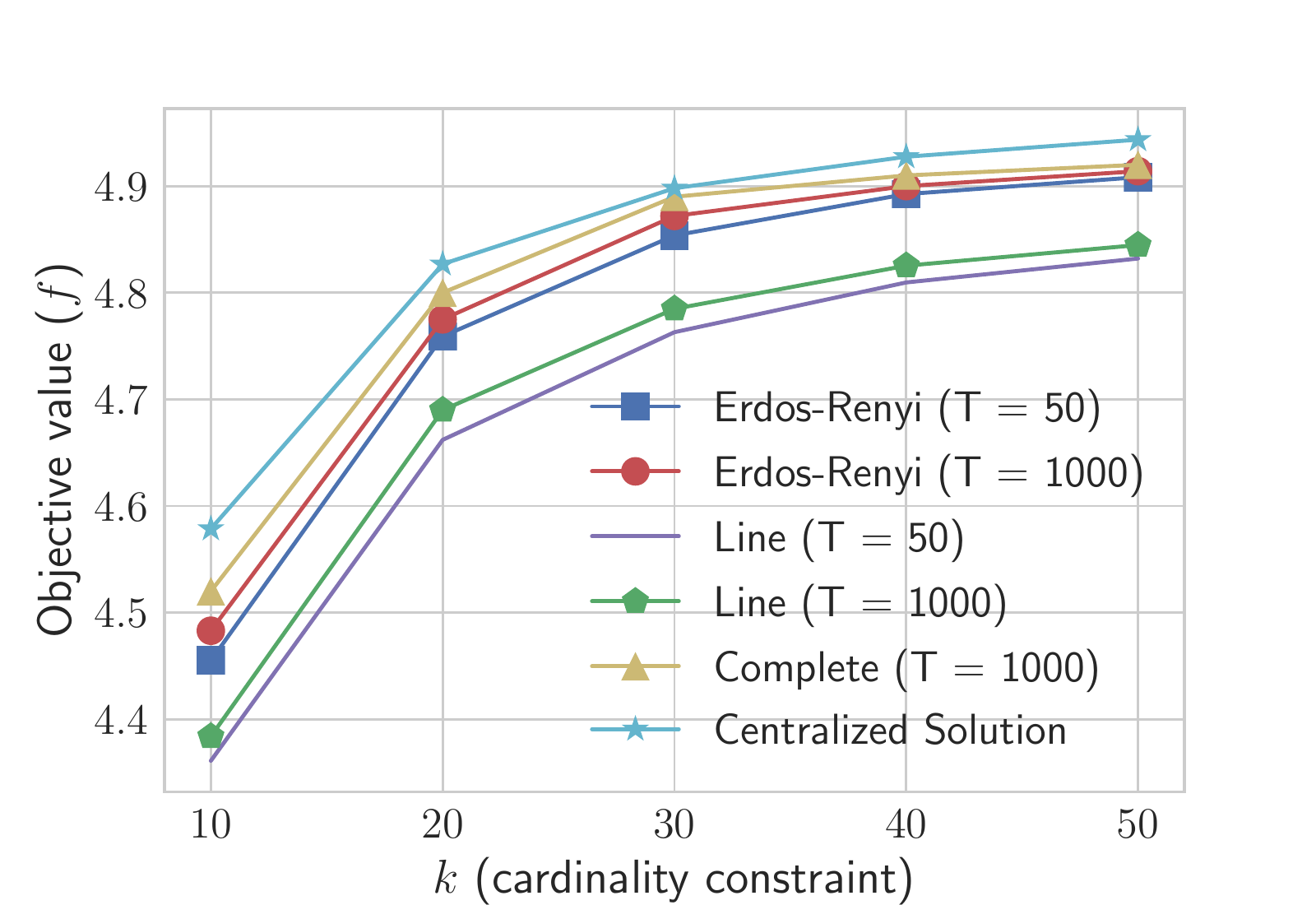}}
\end{center}
%\vspace{-2mm}
\caption{The average objective value is plotted as a function of the cardinality constraint $k$ for different choices of the communication graph as well as number of iterations $T$. Note that ``ER" stands for the Erdos-Reny graph with average degree $5$, ``Line" stands for the line graph and ``Complete" is for the complete graph. We have run Algorithm~\ref{algo_DDCG} for $T = 50$ and $T=1000$. 
}
\label{fig-vals22}
\end{figure} 

Figure~\ref{fig-vals22} depicts the obtained objective value of Discrete DCG (Algorithm~\ref{algo_DDCG}) for the three networks considered above.  More precisely, we plot the value $\frac{1}{n} \sum_{i=1}^n f(\bbx_i^T)$ obtained at the end of Algorithm~\ref{algo_DDCG} as a function of the cardinality constraint $k$.  We also compare these values with the value obtained by the centralized greedy algorithm (i.e. the centralized solution). A few comments are in order. The performance of Algorithm~\ref{algo_DDCG} is close to the centralized solution when the underlying graph is the Erdos-Renyi (with average degree 5) graph or  the complete graphs. This is because for both such graphs consensus is achieved from the early stages of the algorithm. By increasing $T$, we see that the performance becomes closer to the centralized solution. However, when the underlying graph is the line graph, then consensus will not be achieved unless the number of iterations is significantly increased. Consequently, for small number of iterations (e.g., $T \leq 1000$) the performance of the algorithm will not be close to the centralized solution. Indeed, this is not surprising as the line graph is poorly~connected.

\section{Conclusion}\label{sec_conclusion}

In this paper, we proposed the first fully decentralized
optimization method for maximizing a monotone and
continuous DR-submodular function where its components
are available at the different nodes of a connected
graph. We developed Decentralized Continuous
Greedy (DCG) that achieves a $(1-1/e-\eps)$ approximation
guarantee with $\mathcal{O}(1/\eps^2)$ local rounds of communication.
We also showed that our continuous algorithm can be used
to provide the first $(1-1/e)$ tight approximation guarantee
for maximizing a monotone submodular set function subject
to a general matroid constraint in a decentralized fashion. In
particular, we demonstrated that by lifting the local discrete
functions to the continuous domain and using DCG as an interface,
after $\mathcal{O}(1/\eps^3)$ rounds of communication, each node
achieves a tight $(1-1/e-\eps)$ fractional approximate solution.
Such solutions can be efficiently rounded in order to obtain
discrete solutions with the same approximation guarantee.

\acks{The work of A. Mokhtari was partially supported by the DIMACS/Simons Collaboration on Bridging Continuous and Discrete Optimization through NSF grant CCF-1740425. The work of A. Karbasi was supported by DARPA Young Faculty Award (D16AP00046).}

%!TEX root = example_paper.tex
\section{Appendix}

%%%%%%%%%%%%%%%%%%%%%%%%%%%%%%%%%%%
%%%%%%%%%%%%%%%%%%%%%%%%%%%%%%%%%%%
%%%%%%%%%%%%%%%%%%%%%%%%%%%%%%%%%%%
%%%%   S  U  B  S  E  C  T  I  O  N    %%%%%%%%%%%%%%%%
%%%%%%%%%%%%%%%%%%%%%%%%%%%%%%%%%%%
%%%%%%%%%%%%%%%%%%%%%%%%%%%%%%%%%%%

\subsection{Proof of Proposition~\ref{prop_stay}}\label{app:prop_stay}
Define $\bbx_{con}=[\bbx_1;\dots;\bbx_n]\in \reals^{np}$ and $\bbv_{con}=[\bbv_1;\dots;\bbv_n]\in \reals^{np}$ as the concatenation of the local variables and descent directions, respectively. Using these definitions and the update in \eqref{eq:variable_update} we can write
%%%
\begin{align}\label{proof_2_100}
\bbx_{con}^{t+1} = (\bbW\otimes\bbI)\bbx_{con}^{t}+\frac{1}{T} \bbv_{con}^t,
\end{align}
%%%%
where $ \bbW\otimes\bbI\in \reals^{np\times np}$ is the kronecker product of the matrices $\bbW\in\reals^{n\times n}$ and $\bbI\in\reals^{p\times p}$. If we set $\bbx_i^0=\bb0_p$ for all nodes $i$, it follows that $\bbx_{con}^0=\bb0_{np}$. Hence, by applying the update in \eqref{proof_2_100} recursively we obtain that the iterate $\bbx_{con}^t$ is equal to 
%%%
\begin{align}\label{proof_2_200}
\bbx_{con}^{t} =\frac{1}{T} \sum_{s=0}^{t-1}(\bbW\otimes\bbI)^{t-1-s} \bbv_{con}^s.
\end{align}
%%%
%Setting $t=T$ in \eqref{proof_2_200} implies that the concatenation of the local variables after $T$ iterations can be written as 
%%%%%
%%%%
%\begin{align}\label{proof_2_300}
%\bbx_{con}^{T} =\frac{1}{T} \sum_{s=0}^{T-1}(\bbW\otimes\bbI)^{T-1-s} \bbv_{con}^s.
%\end{align}
%%%%

We proceed by showing that if the local blocks of a vector $\bbv_{con}\in\reals^{np}$ belong to the feasible set $\ccalC$, i.e., $\bbv_i\in\ccalC$ for $i=1,\dots,n$, then the local vectors of $\bby_{con}=(\bbW\otimes\bbI)\bbv_{con}\in\reals^{np}$ also in the set $\ccalC$. Note that if the condition $\bby_{con}=(\bbW\otimes\bbI)\bbv_{con}$ holds, then the $i$-th block of $\bby_{con}=[\bby_1;\dots;\bby_n]$ can be written as 
%%%
\begin{align}\label{proof_2_400}
\bby_i=\sum_{j=1}^nw_{ij}\bbv_{j}.
\end{align}
%%%
Since we assume that all $\{\bbv_{j}\}_{j=1}^n$ belong to the set $\ccalC$ and the set $\ccalC$ is convex, the weighted average of these vectors also is in the set $\ccalC$, i.e., $\bby_i\in \ccalC$. This argument indeed holds for all blocks $\bby_i$ and therefore $\bby_i\in \ccalC$ for $i=1,\dots,n$. This argument verifies that if we apply any power of the  matrix $\bbW\otimes\bbI$ to a vector $\bbv_{con}\in\reals^{np}$ whose blocks belong to the set $\ccalC$, then the local components of the output vector also belong to the set $\ccalC$. Therefore, the local components of each of the terms $(\bbW\otimes\bbI)^{t-1-s} \bbv_{con}^s$ in \eqref{proof_2_200} belong to the set $\ccalC$. The fact that $\bbx_i$ which is the $i$-th block of the vector $\bbx_{con}^t$, is the average of $T$ terms that are in the set $\ccalC$  ($\bbx_{con}^t$ is the average of the vectors $(\bbW\otimes\bbI)^{t-1}\bbv_{con}^0, \dots, (\bbW\otimes\bbI)^{0} \bbv_{con}^{t-1} $ with weights $1/T$ and the vector $\bb0_{np}$ with weight $(T-t)/T$), implies that $\bbx_i^t\in \ccalC$.  This result holds for all $i\in\{1,\dots,n\}$ and the proof is complete.

\subsection{Proof of Lemma~\ref{lemma:ar_in_avg_bound}}\label{app:lemma:ar_in_avg_bound}
By averaging both sides of the update in \eqref{eq:variable_update} over the nodes in the network and using the fact $w_{ij}=0$ if $i$ and $j$ are not neighbors we can write
%%%
\begin{align}\label{proof_1_100}
\frac{1}{n}\sum_{i=1}^n \bbx_i^{t+1}
&= \frac{1}{n}\sum_{i=1}^n \sum_{j\in \ccalN_i\cup\{i\}}  w_{ij}\bbx_{j}^{t}+\frac{1}{T} \frac{1}{n}\sum_{i=1}^n\bbv_i^t\nonumber\\
& = \frac{1}{n}\sum_{i=1}^n\sum_{j=1}^n  w_{ij}\bbx_{j}^{t}+\frac{1}{T} \frac{1}{n}\sum_{i=1}^n\bbv_i^t\nonumber\\
& = \frac{1}{n}\sum_{j=1}^n \bbx_{j}^{t} \sum_{i=1}^nw_{ij}+\frac{1}{T} \frac{1}{n}\sum_{i=1}^n\bbv_i^t\nonumber\\
& = \frac{1}{n}\sum_{j=1}^n \bbx_{j}^{t}+\frac{1}{T} \frac{1}{n}\sum_{i=1}^n\bbv_i^t,
\end{align}
%%%%
where the last equality holds since $\bbW^T\bbone_n=\bbone_n$ (i.e. $\bbW$ is a doubly stochastic matrix). By using the definition of the average iterate vector $\bar{\bbx}^t$ and the result in \eqref{proof_1_100} it follows that
%%%
\begin{align}\label{proof_1_200}
\bar{\bbx}^{t+1} = \bar{\bbx}^{t} +\frac{1}{T} \frac{1}{n}\sum_{i=1}^n\bbv_i^t.
\end{align}
%%%
Since $\bbv_i^t$ belongs to the convex set $\ccalC$ its Euclidean norm is bounded by $\|\bbv_i^t\|\leq D$ according to Assumption \ref{ass:bounded_set}. This inequality and the expression in \eqref{proof_1_200} yield 
%%%
\begin{align}
\|\bar{\bbx}^{t+1} - \bar{\bbx}^{t}\| \leq \frac{D}{T},
\end{align}
%%%
and the claim in \eqref{eq:var_in_avg_bound} follows.

%%%%%%%%%%%%%%%%%%%%%%%%%%%%%%%%%%%
%%%%%%%%%%%%%%%%%%%%%%%%%%%%%%%%%%%
%%%%%%%%%%%%%%%%%%%%%%%%%%%%%%%%%%%
%%%%   S  U  B  S  E  C  T  I  O  N    %%%%%%%%%%%%%%%%
%%%%%%%%%%%%%%%%%%%%%%%%%%%%%%%%%%%
%%%%%%%%%%%%%%%%%%%%%%%%%%%%%%%%%%%

\subsection{Proof of Lemma~\ref{lemma:eq:bound_on_dif_from_avg}}\label{app:lemma:eq:bound_on_dif_from_avg}

Recall the definitions $\bbx_{con}=[\bbx_1;\dots;\bbx_n]\in \reals^{np}$ and $\bbv_{con}=[\bbv_1;\dots;\bbv_n]\in \reals^{np}$ for the concatenation of the local variables and descent directions, respectively. These definitions along with the update in \eqref{eq:variable_update} lead to the expression 
%%%
%%%
\begin{align}\label{proof_3_100}
\bbx_{con}^{t} =\frac{1}{T} \sum_{s=0}^{t-1}(\bbW\otimes\bbI)^{t-1-s} \bbv_{con}^s.
\end{align}
%%%
If we premultiply both sides of \eqref{proof_3_100} by the matrix $(\frac{\bbone_n\bbone_n^{\dag}}{n}\otimes\bbI) $ which is the kronecker product of the matrices $(1/n)(\bbone_n\bbone_n^{\dag})\in\reals^{n\times n}$ and $\bbI\in\reals^{p\times p}$ we obtain
%%%
%%%
\begin{align}\label{proof_3_200}
\left(\frac{\bbone_n\bbone_n^{\dag}}{n}\otimes \bbI\right) \bbx_{con}^t =\frac{1}{T} \sum_{s=0}^{t-1}\left(\left(\frac{\bbone_n\bbone_n^{\dag}}{n}\bbW^{t-1-s}\right)\otimes\bbI\right) \bbv_{con}^s.
\end{align}
%%%
The left hand side of \eqref{proof_3_200} can be simplified to
%%%
\begin{align}\label{proof_3_300}
\left(\frac{\bbone_n\bbone_n^{\dag}}{n}\otimes \bbI\right) \bbx_{con}^t=\bar{\bbx}_{con}^t ,
\end{align}
%%%
where $\bar{\bbx}_{con}^t=[\bar{\bbx}^t;\dots;\bar{\bbx}^t]$ is the concatenation of $n$ copies of the average vector $\bar{\bbx}^t$. Using the equality in \eqref{proof_3_300} and the simplification $\bbone_n\bbone_n^{\dag}\bbW=\bbone_n\bbone_n^{\dag}$, we can rewrite  \eqref{proof_3_200} as
%%%
\begin{align}\label{proof_3_400}
\bar{\bbx}_{con}^t= \frac{1}{T} \sum_{s=0}^{t-1}\left(\frac{\bbone_n\bbone_n^{\dag}}{n}\otimes\bbI\right) \bbv_{con}^s.
\end{align}
%%%
Using the expressions in \eqref{proof_3_100} and \eqref{proof_3_400} we can derive an upper bound on the difference $ \|\bbx_{con}^{t}-\bar{\bbx}_{con}^t\|$ as
\begin{align}\label{proof_3_500}
 \|\bbx_{con}^{t}-\bar{\bbx}_{con}^t\|
& =\frac{1}{T} \left\| \sum_{s=0}^{t-1}\left(\left[ \bbW^{t-1-s}-\frac{\bbone_n\bbone_n^{\dag}}{n} \right]\otimes\bbI\right) \bbv_{con}^s \right\|\nonumber\\
& \leq\frac{1}{T} \sum_{s=0}^{t-1}\left\| \bbW^{t-1-s}-\frac{\bbone_n\bbone_n^{\dag}}{n}  \right\|\|\bbv_{con}^s\|\nonumber\\
& \leq\frac{\sqrt{n}D}{T} \sum_{s=0}^{t-1}\left\| \bbW^{t-1-s}-\frac{\bbone_n\bbone_n^{\dag}}{n}  \right\|,
\end{align}
%%%
where the first inequality follows from the Cauchy-Schwarz inequality and 
the fact that the norm of a matrix does not change if we kronecker it by the identity matrix, the second inequality holds since $\|\bbv_i^t\|\leq D$ and therefore $\|\bbv_{con}^t\|\leq \sqrt{n}D$. Note that the eigenvectors of the matrices $\bbW$ and $\bbW^{t-s-1}$ are the same for all $s=0,\dots,t-1$. Therefore, the largest eigenvalue of  $\bbW^{t-s-1}$ is 1 with eigenvector $\bbone_n$ and its second largest magnitude of the eigenvalues is $\beta^{t-1-s}$, where $\beta$ is the second largest magnitude of the eigenvalues of $\bbW$. Also, note that since $\bbW^{t-1-s}$ has $\bbone_n$ as one of its eigenvectors, then all the other eigenvectors of $\bbW$ are orthogonal to $\bbone_n$.  Hence, we can bound the norm $\|\bbW^{t-1-s}-({\bbone_n\bbone_n^{\dag}})/(n)\|$ by $\beta^{t-1-s}$. Applying this substitution into the right hand side of \eqref{proof_3_500} yields 
\begin{align}\label{proof_3_600}
 \|\bbx_{con}^{t}-\bar{\bbx}_{con}^t\|
 \leq\frac{\sqrt{n}D}{T} \sum_{s=0}^{t-1} \beta^{t-1-s} \leq\frac{\sqrt{n}D}{T(1-\beta)}.
\end{align}
%%%
Since $ \|\bbx_{con}^{t}-\bar{\bbx}_{con}^t\|^2=\sum_{i=1}^n \left\|\bbx_i^t-\bar{\bbx}^t\right\|^2$, the claim in \eqref{eq:bound_on_dif_from_avg} follows.

%%%%%%%%%%%%%%%%%%%%%%%%%%%%%%%%%%%
%%%%%%%%%%%%%%%%%%%%%%%%%%%%%%%%%%%
%%%%%%%%%%%%%%%%%%%%%%%%%%%%%%%%%%%
%%%%   S  U  B  S  E  C  T  I  O  N    %%%%%%%%%%%%%%%%
%%%%%%%%%%%%%%%%%%%%%%%%%%%%%%%%%%%
%%%%%%%%%%%%%%%%%%%%%%%%%%%%%%%%%%%
\subsection{Proof of Lemma \ref{lemma:bound_on_gradient_consensus_error}}\label{proof:lemma:bound_on_gradient_consensus_error}

Recall the definition of the vector $\bbx_{con}=[\bbx_1;\dots;\bbx_n]\in \reals^{np}$ as the concatenation of the local variables, and define $\bbd_{con}=[\bbd_1;\dots;\bbd_n]\in \reals^{np}$ as the concatenation of the local approximate gradients. Further, consider the function $F_{con}:\ccalX^n\to\reals$ which is defined as $F_{con}(\bbx_{con})=F_{con}(\bbx_1,\dots,\bbx_n):=\sum_{i=1}^n F_i(\bbx_i)$. According to these definitions and the update in \eqref{eq:gradient_update}, we can show that 
%%%%
\begin{align}\label{eq:proof_10_1}
\bbd_{con}^t = (1-\alpha) (\bbW\otimes\bbI) \bbd_{con}^{t-1}+\alpha \nabla F_{con}(\bbx_{con}^t),
\end{align}
%%%%
where $ \bbW\otimes\bbI\in \reals^{np\times np}$ is the kronecker product of the matrices $\bbW\in\reals^{n\times n}$ and $\bbI\in\reals^{p\times p}$. Considering the initialization $\bbd_{con}^0=\bb0_p$, applying the update in \eqref{eq:proof_10_1} recursively from step $1$ to $t$ leads to 
%%%%
\begin{align}\label{eq:proof_10_2}
\bbd_{con}^t=\alpha \sum_{s=1}^{t} (((1-\alpha)\bbW)^{t-s}\otimes \bbI) \nabla F_{con}(\bbx_{con}^s).
\end{align}
%%%%
If we multiply both sides of \eqref{eq:proof_10_2} from left by the matrix $(\frac{\bbone_n\bbone_n^{\dag}}{n}\otimes\bbI) \in \reals^{np\times np}$ and use the properties of the weight matrix $\bbW$, i.e., $\bbone_n^{\dag}\bbW^{t-s}=\bbone_n^{\dag}$, we obtain that
%%%%
\begin{align}\label{eq:proof_10_3}
\bar{\bbd}_{con}^t=\alpha \sum_{s=1}^{t} (1-\alpha)^{t-s}\left(\frac{\bbone_n\bbone_n^{\dag}}{n}\otimes \bbI\right) \nabla F_{con}(\bbx_{con}^s),
\end{align}
%%%%
where $\bar{\bbd}_{con}^t=[\bar{\bbd}^t;\dots;\bar{\bbd}^t]$ is the concatenation of $n$ copies of the average vector $\bar{\bbd}^t$. Hence, the difference $ \|\bbd_{con}^t -\bar{\bbd}_{con}^t\|$ can be upper bounded by
%%%%
\begin{align}\label{eq:proof_10_4}
\|\bbd_{con}^t -\bar{\bbd}_{con}^t\|  &= \alpha \left\|  \sum_{s=1}^{t} (1-\alpha)^{t-s}(\bbW^{t-s}\otimes \bbI) \nabla F_{con}(\bbx_{con}^s)- \sum_{s=1}^{t} (1-\alpha)^{t-s}\left[\frac{\bbone_n\bbone_n^{\dag}}{n}\otimes \bbI\right] \nabla F_{con}(\bbx_{con}^s)  \right\|  \nonumber\\
&= \alpha \left\|  \sum_{s=1}^{t}(1-\alpha)^{t-s} \left[(\bbW^{t-s}-\frac{\bbone_n\bbone_n^{\dag}}{n})\otimes \bbI \right] \nabla F_{con}(\bbx_{con}^s)  \right\|  \nonumber\\
&\leq   \alpha \sqrt{n}G \sum_{s=1}^{t}(1-\alpha)^{t-s} \beta^{t-s}  \nonumber\\
&\leq   \frac{\alpha \sqrt{n} G}{1-\beta(1-\alpha)},
\end{align}
%%%%
where the first equality is implied by replacing $\bbd_{con}^t$ and $\bar{\bbd}_{con}^t$ with the expressions in \eqref{eq:proof_10_2} and \eqref{eq:proof_10_3}, respectively, the second equality is achieved by regrouping the terms,
 the first inequality holds since $\|\nabla F_i(x_i^s)\|\leq G$ and $\|\bbW^{t-s-1}-(\bbone_n\bbone_n^{\dag})/n\|\leq \beta^{t-s-1}$, and finally the last inequality is valid since $\sum_{s=1}^{t}((1-\alpha)\beta)^{t-s}\leq \frac{1}{1-(\beta(1-\alpha))}$. Now considering the result in \eqref{eq:proof_10_4} and the expression $ \|\bbd_{con}^t -\bar{\bbd}_{con}^t\|^2=\sum_{i=1}^n \|\bbd_i^t-\bar{\bbd}^t\|^2$, the claim in \eqref{eq:bound_on_gradient_consensus_error} follows.

%%%%%%%%%%%%%%%%%%%%%%%%%%%%%%%%%%%
%%%%%%%%%%%%%%%%%%%%%%%%%%%%%%%%%%%
%%%%%%%%%%%%%%%%%%%%%%%%%%%%%%%%%%%
%%%%   S  U  B  S  E  C  T  I  O  N    %%%%%%%%%%%%%%%%
%%%%%%%%%%%%%%%%%%%%%%%%%%%%%%%%%%%
%%%%%%%%%%%%%%%%%%%%%%%%%%%%%%%%%%%
\subsection{Proof of Lemma \ref{lemma:bound_on_gradient_error}}\label{proof:lemma:bound_on_gradient_error}

Considering the update in \eqref{eq:gradient_update}, we can write the sum of local ascent directions $\bbd_i^t$ at step $t$ as
%%%%
\begin{align}\label{app_proof_100}
\sum_{i=1}^n \bbd_i^t 
&= (1-\alpha) \sum_{i=1}^n \sum_{j=1}^n  w_{ij}\bbd_{j}^{t-1}+\alpha \sum_{i=1}^n \nabla F_i(\bbx_i^t) \nonumber\\
& = (1-\alpha)  \sum_{j=1}^n \bbd_{j}^{t-1}  \sum_{i=1}^n w_{ij}+\alpha \sum_{i=1}^n \nabla F_i(\bbx_i^t) \nonumber\\
&= (1-\alpha)  \sum_{j=1}^n \bbd_{j}^{t-1} +\alpha \sum_{i=1}^n \nabla F_i(\bbx_i^t),
\end{align}
%%%%
where the last equality holds since $\sum_{i=1}^n w_{ij}=1$ which is the consequence of $\bbW^\dag \bbone_n=\bbone_n$. Now, we use the expression in \eqref{app_proof_100} to bound the difference $\|\sum_{i=1}^n \bbd_i^t  - \sum_{i=1}^n\nabla F_i(\bar{\bbx}^t) \|$. Hence,
%%%
\begin{align}\label{app_proof_200}
&\left\|\sum_{i=1}^n \bbd_i^t  - \sum_{i=1}^n\nabla F_i(\bar{\bbx}^t) \right\|\nonumber\\
&=\left\| (1-\alpha)  \sum_{j=1}^n \bbd_{j}^{t-1}  +\alpha \sum_{i=1}^n \nabla F_i(\bbx_i^t)  - \sum_{i=1}^n\nabla F_i(\bar{\bbx}^t) \right\| \nonumber\\
&=\left\| (1-\alpha)  \sum_{j=1}^n \bbd_{j}^{t-1}  -(1-\alpha ) \sum_{i=1}^n \nabla F_i(\bar{\bbx}^{t-1})  + (1-\alpha)\sum_{i=1}^n\nabla F_i(\bar{\bbx}^{t-1})+\alpha \sum_{i=1}^n \nabla F_i(\bbx_i^t)  - \sum_{i=1}^n\nabla F_i(\bar{\bbx}^t)  \right\| \nonumber\\
&=\Bigg\| (1-\alpha) \!\left[ \sum_{j=1}^n \bbd_{j}^{t-1}-\sum_{i=1}^n\nabla F_i(\bar{\bbx}^{t-1})  \right]\!
 +(1-\alpha )\! \left[\sum_{i=1}^n \nabla F_i(\bar{\bbx}^{t-1})  - \sum_{i=1}^n \nabla F_i(\bar{\bbx}^{t}) \right]\!
\nonumber\\ 
&\qquad  +\alpha\! \left[\sum_{i=1}^n   \nabla F_i(\bbx_i^t)  - \nabla F_i(\bar{\bbx}^{t})   \right]\!
 \Bigg\| \nonumber\\
 &\leq (1-\alpha) \left\| \sum_{j=1}^n \bbd_{j}^{t-1}-\sum_{i=1}^n\nabla F_i(\bar{\bbx}^{t-1})  \right\|
 +(1-\alpha ) \left\|\sum_{i=1}^n \nabla F_i(\bar{\bbx}^{t-1})  - \sum_{i=1}^n \nabla F_i(\bar{\bbx}^{t}) \right\|
 \nonumber\\
&\qquad  +\alpha \left\|\sum_{i=1}^n   \nabla F_i(\bbx_i^t)  - \nabla F_i(\bar{\bbx}^{t})  \right\| .
\end{align}
The first equality is the outcome of replacing $\sum_{i=1}^n \bbd_i^t $ by the expression in \eqref{app_proof_100}, the second equality is obtained by adding and subtracting $(1-\alpha ) \sum_{i=1}^n \nabla F_i(\bar{\bbx}^{t-1})$, in the third equality we regroup the terms, and the inequality follows from applying the triangle inequality twice. Applying the Cauchy--Schwarz inequality to the second and third summands in \eqref{app_proof_200} and using the Lipschitz continuity of the gradients lead to 
%%%
\begin{align}\label{app_proof_300}
&\left\|\sum_{i=1}^n \bbd_i^t  - \sum_{i=1}^n\nabla F_i(\bar{\bbx}^t) \right\| \nonumber\\
&\leq (1-\alpha) \left\|  \sum_{j=1}^n \bbd_{j}^{t-1}-\sum_{i=1}^n\nabla F_i(\bar{\bbx}^{t-1}) \right\| 
+(1-\alpha ) L \sum_{i=1}^n \left\|\bar{\bbx}^{t-1} - \bar{\bbx}^{t}\right\| 
+ \alpha L \sum_{i=1}^n \|  \bbx_i^t  - \bar{\bbx}^{t}\|.
\end{align}
%%%
According to the result in Lemma \ref{lemma:ar_in_avg_bound}, we can bound the $\sum_{i=1}^n\|\bar{\bbx}^{t+1} - \bar{\bbx}^{t}\| $ by $ {nD}/{T}$. Further, the result in Lemma \ref{lemma:eq:bound_on_dif_from_avg} shows that $(\sum_{i=1}^n \|\bbx_i^t-\bar{\bbx}^t\|^2)^{1/2} \leq  \frac{\sqrt{n}D}{T(1-\beta)}.$ Since by the Cauchy--Swartz inequality it holds that $(\sum_{i=1}^n \|\bbx_i^t-\bar{\bbx}^t\|^2)^{1/2}\geq \frac{1}{\sqrt{n}}\sum_{i=1}^n \|\bbx_i^t-\bar{\bbx}^t\|$, it follows that $\sum_{i=1}^n \|\bbx_i^t-\bar{\bbx}^t\|\leq (nD)/(T(1-\beta))$. Applying these substitutions into \eqref{app_proof_300} yields
%%%
\begin{align}\label{app_proof_400}
\left\|\sum_{i=1}^n \bbd_i^t  - \sum_{i=1}^n\nabla F_i(\bar{\bbx}^t) \right\|
\leq (1-\alpha) \left\|  \sum_{j=1}^n \bbd_{j}^{t-1}-\sum_{i=1}^n\nabla F_i(\bar{\bbx}^{t-1}) \right\| 
+\frac{(1-\alpha )LnD}{T} 
+ \frac{\alpha LnD}{T(1-\beta)} .
\end{align}
By multiplying both of sides of \eqref{app_proof_400} by $1/n$ and applying the resulted inequality recessively for $t$ steps we obtain 
%%%
\begin{align}\label{app_proof_500}
&\left\|\frac{1}{n}\sum_{i=1}^n \bbd_i^t  - \frac{1}{n}\sum_{i=1}^n\nabla F_i(\bar{\bbx}^t) \right\|\nonumber\\
&\leq (1-\alpha)^t \left\| \frac{1}{n} \sum_{j=1}^n \bbd_{j}^{0}-\frac{1}{n}\sum_{i=1}^n\nabla F_i(\bar{\bbx}^{0}) \right\| 
+\left(\frac{(1-\alpha )LD}{T} 
+ \frac{\alpha LD}{T(1-\beta)} \right)\sum_{s=0}^{t-1} (1-\alpha)^s\nonumber\\
&\leq (1-\alpha)^t \frac{1}{n}\sum_{i=1}^n\left\|\nabla F_i(\bar{\bbx}^{0}) \right\| 
+\frac{(1-\alpha)LD}{\alpha T} 
+ \frac{ LD}{T(1-\beta)} \nonumber\\
&\leq (1-\alpha)^t G
+\frac{(1-\alpha)LD}{\alpha T} 
+ \frac{ LD}{T(1-\beta)} ,
\end{align}
where the second inequality holds since $ \sum_{j=1}^n \bbd_{j}^{0}=\bb0_p$ and $\sum_{s=0}^{t-1} (1-\alpha)^s\leq 1/\alpha$, and the last inequality follows from Assumption \ref{ass:smoothness2}. 

%%%%%%%%%%%%%%%%%%%%%%%%%%%%%%%%%%%
%%%%%%%%%%%%%%%%%%%%%%%%%%%%%%%%%%%
%%%%%%%%%%%%%%%%%%%%%%%%%%%%%%%%%%%
%%%%   S  U  B  S  E  C  T  I  O  N    %%%%%%%%%%%%%%%%
%%%%%%%%%%%%%%%%%%%%%%%%%%%%%%%%%%%
%%%%%%%%%%%%%%%%%%%%%%%%%%%%%%%%%%%

\subsection{Proof of Theorem \ref{theorem:main_theorem}}\label{proof:theorem:main_theorem}

Recall the definition of  $\bar{\bbx}^{t}=\frac{1}{n}\sum_{i=1}^n \bbx_i^t$ as the average of local variables at step $t$. Since the gradients of the global objective function are $L$-Lipschitz we can write 
%%%%
\begin{align}\label{final_thm_proof_100}
\frac{1}{n}\sum_{i=1}^n F_i(\bar{\bbx}^{t+1})-\frac{1}{n}\sum_{i=1}^n F_i(\bar{\bbx}^{t}) 
&\geq \frac{1}{n}\langle \sum_{i=1}^n \nabla F_i(\bar{\bbx}^{t}), \bar{\bbx}^{t+1}-\bar{\bbx}^{t} \rangle -\frac{L}{2}\| \bar{\bbx}^{t+1}-\bar{\bbx}^{t}\|^2 \nonumber\\
&= \frac{1}{T}\langle \frac{1}{n}\sum_{i=1}^n \nabla F_i(\bar{\bbx}^{t}),  \frac{1}{n}\sum_{i=1}^n\bbv_i^t \rangle -\frac{L}{2T^2}\left\|\frac{1}{n}\sum_{i=1}^n\bbv_i^t \right\|^2,
\end{align}
%%%
where the equality holds due to the expression in \eqref{proof_1_200}. Note that the term $\|(1/n)\sum_{i=1}^n\bbv_i^t \|^2$ can be upper bounded by $D^2$ according to Assumption \ref{ass:bounded_set}, since $(1/n)\sum_{i=1}^n\bbv_i^t \in \ccalC$ . Apply this substition into \eqref{final_thm_proof_100} and add and subtract $(1/nT)\sum_{i=1}^n \bbd_i^t$ to obtain
%%%%
\begin{align}\label{final_thm_proof_200}
&\frac{1}{n}\sum_{i=1}^n F_i(\bar{\bbx}^{t+1})-\frac{1}{n}\sum_{i=1}^n F_i(\bar{\bbx}^{t}) \nonumber\\
&\geq \frac{1}{T}\langle\frac{1}{n} \sum_{i=1}^n  \bbd_i^t,  \frac{1}{n}\sum_{i=1}^n\bbv_i^t \rangle +\frac{1}{T}\langle \frac{1}{n}\sum_{i=1}^n \nabla F_i(\bar{\bbx}^{t})-\frac{1}{n}\sum_{i=1}^n \bbd_i^t,  \frac{1}{n}\sum_{i=1}^n\bbv_i^t \rangle
-\frac{LD^2}{2T^2}.
\end{align}
%%%%
Now by rewriting the inner product $\langle \sum_{i=1}^n  \bbd_i^t,\sum_{i=1}^n\bbv_i^t \rangle $ as $\sum_{i=1}^n\sum_{j=1}^n \langle   \bbd_i^t, \bbv_j^t \rangle= \sum_{j=1}^n \langle  \sum_{i=1}^n \bbd_i^t, \bbv_j^t \rangle$, we can rewrite the right hand side of \eqref{final_thm_proof_200} as
%%%
\begin{align}\label{final_thm_proof_300}
&\frac{1}{n}\sum_{i=1}^n F_i(\bar{\bbx}^{t+1})-\frac{1}{n}\sum_{i=1}^n F_i(\bar{\bbx}^{t}) \nonumber\\
&\geq \frac{1}{n^2T} \sum_{j=1}^n  \langle  \sum_{i=1}^n\bbd_i^t,  \bbv_j^t \rangle 
+\frac{1}{T}\langle\frac{1}{n} \sum_{i=1}^n \nabla F_i(\bar{\bbx}^{t})-\frac{1}{n}\sum_{i=1}^n \bbd_i^t,  \frac{1}{n}\sum_{i=1}^n\bbv_i^t \rangle
-\frac{LD^2}{2T^2}\nonumber\\
&= \frac{1}{nT} \sum_{j=1}^n  \langle \bbd_j^t,  \bbv_j^t \rangle 
+\frac{1}{nT} \sum_{j=1}^n  \langle (\frac{1}{n} \sum_{i=1}^n\bbd_i^t-\bbd_j^t),  \bbv_j^t \rangle 
+\frac{1}{T}\langle \sum_{i=1}^n \frac{1}{n}\nabla F_i(\bar{\bbx}^{t})-\frac{1}{n}\sum_{i=1}^n \bbd_i^t,  \frac{1}{n}\sum_{i=1}^n\bbv_i^t \rangle
-\frac{LD^2}{2T^2}.
\end{align} 
%%%%
Note that in the last step we added and and subtracted $ (1/nT) \sum_{j=1}^n  \langle \bbd_j^t,  \bbv_j^t \rangle$. Now according to the update in \eqref{eq:descent_update} we can write, $\langle \bbd_j^t,  \bbv_j^t \rangle =  \max_{\bbv\in \ccalC}\ \langle \bbd_j^t, \bbv\rangle \geq   \langle \bbd_j^t, \bbx^*\rangle$. Hence, we can replace $\langle \bbd_j^t,  \bbv_j^t \rangle$ by its lower bound $ \langle \bbd_j^t, \bbx^*\rangle$ to obtain
%%%
\begin{align} \label{final_thm_proof_400}
&\frac{1}{n}\sum_{i=1}^n F_i(\bar{\bbx}^{t+1})-\frac{1}{n}\sum_{i=1}^n F_i(\bar{\bbx}^{t})  \nonumber\\
&\geq \frac{1}{nT} \sum_{j=1}^n  \langle \bbd_j^t,  \bbx^* \rangle 
+\frac{1}{nT} \sum_{j=1}^n  \langle ( \frac{1}{n}\sum_{i=1}^n\bbd_i^t-\bbd_j^t),  \bbv_j^t \rangle 
+\frac{1}{T}\langle \frac{1}{n}\sum_{i=1}^n \nabla F_i(\bar{\bbx}^{t})-\frac{1}{n}\sum_{i=1}^n \bbd_i^t,  \frac{1}{n}\sum_{i=1}^n\bbv_i^t \rangle
-\frac{LD^2}{2T^2}.
\end{align}
%%%%
Adding and subtracting $\frac{1}{n^2T} \sum_{j=1}^n  \langle \sum_{i=1}^n \bbd_i^t,  \bbx^* \rangle $ and regrouping the terms lead to
%%%%
\begin{align} \label{final_thm_proof_500}
\frac{1}{n}\sum_{i=1}^n F_i(\bar{\bbx}^{t+1})&-\frac{1}{n}\sum_{i=1}^n F_i(\bar{\bbx}^{t})  
\geq \frac{1}{n^2T} \sum_{j=1}^n  \langle \sum_{i=1}^n \bbd_i^t,  \bbx^* \rangle 
+\frac{1}{nT} \sum_{j=1}^n  \langle \bbd_j^t -\frac{1}{n}\sum_{i=1}^n \bbd_i^t,  \bbx^* \rangle 
\nonumber\\
&\qquad  +\frac{1}{nT} \sum_{j=1}^n  \langle  \frac{1}{n}\sum_{i=1}^n\bbd_i^t-\bbd_j^t,  \bbv_j^t \rangle +\frac{1}{T}\langle \frac{1}{n}\sum_{i=1}^n \nabla F_i(\bar{\bbx}^{t})-\frac{1}{n}\sum_{i=1}^n \bbd_i^t,  \frac{1}{n}\sum_{i=1}^n\bbv_i^t \rangle
-\frac{LD^2}{2T^2}.
\end{align}
%%%%
Further add and subtract the expression $\frac{1}{n^2T} \sum_{j=1}^n  \langle  \sum_{i=1}^n \nabla F_i(\bar{\bbx}^t),  \bbx^* \rangle $ and combine the terms to obtain
%%%%
\begin{align} \label{final_thm_proof_600}
&\frac{1}{n}\sum_{i=1}^n F_i(\bar{\bbx}^{t+1})-\frac{1}{n}\sum_{i=1}^n F_i(\bar{\bbx}^{t})  \nonumber\\
&\geq \frac{1}{n^2T} \sum_{j=1}^n  \langle  \sum_{i=1}^n \nabla F_i(\bar{\bbx}^t),  \bbx^* \rangle 
+\frac{1}{nT} \sum_{j=1}^n  \langle ( \frac{1}{n}\sum_{i=1}^n \bbd_i^t - \frac{1}{n} \sum_{i=1}^n \nabla F_i(\bar{\bbx}^t),  \bbx^* \rangle 
+\frac{1}{nT} \sum_{j=1}^n  \langle (\bbd_j^t - \frac{1}{n}\sum_{i=1}^n \bbd_i^t,  \bbx^* \rangle 
\nonumber\\
&\qquad 
+\frac{1}{nT} \sum_{j=1}^n  \langle ( \frac{1}{n} \sum_{i=1}^n\bbd_i^t-\bbd_j^t),  \bbv_j^t \rangle 
+\frac{1}{T}\langle \sum_{i=1}^n \frac{1}{n} \nabla F_i(\bar{\bbx}^{t})- \frac{1}{n}\sum_{i=1}^n \bbd_i^t,  \frac{1}{n}\sum_{i=1}^n\bbv_i^t \rangle
-\frac{LD^2}{2T^2}\nonumber\\
&= \frac{1}{nT}   \langle  \sum_{i=1}^n \nabla F_i(\bar{\bbx}^t),  \bbx^* \rangle 
+\frac{1}{nT}  \langle \sum_{i=1}^n \bbd_i^t - \sum_{i=1}^n \nabla F_i(\bar{\bbx}^t),  \bbx^*-\frac{1}{n}\sum_{i=1}^n\bbv_i^t \rangle \nonumber\\
&\qquad +\frac{1}{nT} \sum_{j=1}^n  \langle ( \frac{1}{n}\sum_{i=1}^n\bbd_i^t-\bbd_j^t),  \bbv_j^t-\bbx^* \rangle 
-\frac{LD^2}{2T^2}.
\end{align}
%%%
The monotonicity of the average function $(1/n)\sum_{i=1}^n F_i(\bbx)$ and its concavity along positive directions imply that $\langle  (1/n) \sum_{i=1}^n \nabla F_i(\bar{\bbx}^t),  \bbx^* \rangle \geq (1/n) \sum_{i=1}^n F_i(\bbx^*)- (1/n) \sum_{i=1}^n F_i(\bar{\bbx}^t)$. By applying this substitution into \eqref{final_thm_proof_600} and using the Cauchy-Schwarz inequality we obtain 
%%%%
\begin{align} \label{final_thm_proof_700}
\frac{1}{n}\sum_{i=1}^n F_i(\bar{\bbx}^{t+1})-\frac{1}{n}\sum_{i=1}^n F_i(\bar{\bbx}^{t})  
&\geq \frac{1}{nT} \left[ \sum_{i=1}^n F_i(\bbx^*)- \sum_{i=1}^n F_i(\bar{\bbx}^t)\right]
\!-\!\frac{1}{nT} \sum_{j=1}^n  \left\| \frac{1}{n}\sum_{i=1}^n\bbd_i^t-\bbd_j^t\right\|  \|  \bbv_j^t-\bbx^* \|
\nonumber\\
&\qquad \quad 
-\frac{1}{nT}   \left\|\sum_{i=1}^n \bbd_i^t - \sum_{i=1}^n \nabla F_i(\bar{\bbx}^t) \right\|   \left\|  \bbx^*-\frac{1}{n}\sum_{i=1}^n\bbv_i^t \right\|-\frac{LD^2}{2T^2}.
\end{align}
%%%%
Now we proceed to derive lower bounds for the negative terms on the right hand side of \eqref{final_thm_proof_700}. Note that all $\bbv_i^t$ for $i=1,\dots,n$ belong to the convex set $\ccalC$ and therefore the average vector $\frac{1}{n}\sum_{i=1}^n\bbv_i^t $ is also in the set. Hence, we can bound the difference  $\|  \bbx^*-\frac{1}{n}\sum_{i=1}^n\bbv_i^t \|$ by $D$ according to Assumption \ref{ass:bounded_set}. Indeed, the norm $\|  \bbv_j^t-\bbx^* \|$ is also upper bounded by $D$ and hence we can write
%%%%
\begin{align} \label{final_thm_proof_800}
&\frac{1}{n}\sum_{i=1}^n F_i(\bar{\bbx}^{t+1})-\frac{1}{n}\sum_{i=1}^n F_i(\bar{\bbx}^{t})  \nonumber\\
&\geq \frac{1}{nT} \left[ \sum_{i=1}^n F_i(\bbx^*)- \sum_{i=1}^n F_i(\bar{\bbx}^t)\right]
-\frac{D}{nT}   \left\|\sum_{i=1}^n \bbd_i^t - \sum_{i=1}^n \nabla F_i(\bar{\bbx}^t) \right\|   
-\frac{D}{nT} \sum_{j=1}^n  \left\| \frac{1}{n}\sum_{i=1}^n\bbd_i^t-\bbd_j^t\right\|  
-\frac{LD^2}{2T^2}.
\end{align}
%%%%
The result in Lemma \ref{lemma:bound_on_gradient_consensus_error} implies that 
$(\sum_{i=1}^n\|\bbd_i^t-\bar{\bbd}^t\|^2)^{1/2} \leq \frac{\alpha \sqrt{n} G}{1-\beta(1-\alpha)}$. Note that based on the Cauchy--Swartz inequality it holds that $(\sum_{i=1}^n\|\bbd_i^t-\bar{\bbd}^t\|^2)^{1/2}\geq \frac{1}{\sqrt{n}}\sum_{i=1}^n\|\bbd_i^t-\bar{\bbd}^t\|$, and hence, $\sum_{i=1}^n\|\bbd_i^t-\bar{\bbd}^t\|\leq  \frac{\alpha n G}{1-\beta(1-\alpha)}$. Using this result and recalling the definition $\bar{\bbd}^t:=(1/n)\sum_{i=1}^n\bbd_i^t$, we obtain that 
%%%
\begin{equation}\label{final_thm_proof_900}
\frac{1}{n} \sum_{j=1}^n  \left\| \frac{1}{n}\sum_{i=1}^n\bbd_i^t-\bbd_j^t\right\|  \leq \frac{\alpha G}{1-\beta(1-\alpha)}. 
\end{equation}
%%%
Replace the term $\frac{1}{n} \sum_{j=1}^n  \left\| \frac{1}{n}\sum_{i=1}^n\bbd_i^t-\bbd_j^t\right\|$ in \eqref{final_thm_proof_800} by its upper bound in \eqref{final_thm_proof_900} and use the result in Lemma \ref{lemma:bound_on_gradient_error} to replace $\frac{1}{n}\|\sum_{i=1}^n \bbd_i^t  - \sum_{i=1}^n\nabla F_i(\bar{\bbx}^t) \|$ by its upper bound in \eqref{eq:bound_on_gradient_error}. Applying these substitutions yields
%%%
\begin{align} \label{final_thm_proof_1000}
\frac{1}{n}\sum_{i=1}^n F_i(\bar{\bbx}^{t+1})-\frac{1}{n}\sum_{i=1}^n F_i(\bar{\bbx}^{t}) 
&\geq \frac{1}{T} \left[ \frac{1}{n}\sum_{i=1}^n F_i(\bbx^*)- \frac{1}{n}\sum_{i=1}^n F_i(\bar{\bbx}^t)\right]
-   \frac{(1-\alpha)^t GD}{T}
-\frac{(1-\alpha)LD^2}{\alpha T^2}  \nonumber\\
& \qquad - \frac{ LD^2}{(1-\beta)T^2} 
- \frac{\alpha GD}{(1-\beta(1-\alpha))T}
-\frac{LD^2}{2T^2}.
\end{align}
%%%
Set $\alpha=1/\sqrt{T}$ and regroup the terms to obtain
%%%
\begin{align} \label{final_thm_proof_1100}
 \frac{1}{n}\sum_{i=1}^n F_i(\bbx^*)-\frac{1}{n}\sum_{i=1}^n F_i(\bar{\bbx}^{t+1})
&\leq\left(1- \frac{1}{T}\right) \left[ \frac{1}{n}\sum_{i=1}^n F_i(\bbx^*)- \frac{1}{n}\sum_{i=1}^n F_i(\bar{\bbx}^t)\right]
+ \frac{(1-(1/\sqrt{T}))^t GD}{T}
 \nonumber\\
&\qquad +\frac{LD^2}{ T^{3/2}} + \frac{ LD^2}{(1-\beta)T^2} 
+ \frac{ GD}{(1-\beta)T^{3/2}}
+\frac{LD^2}{2T^2}.
\end{align}
%%%%
By applying the inequality in \eqref{final_thm_proof_1100} recursively for $t =0,\dots,T-1$ we obtain
%%%%
\begin{align} \label{final_thm_proof_1200}
 \frac{1}{n}\sum_{i=1}^n F_i(\bbx^*)-\frac{1}{n}\sum_{i=1}^n F_i(\bar{\bbx}^{T}) 
&\leq\left(1- \frac{1}{T}\right)^T \left[ \frac{1}{n}\sum_{i=1}^n F_i(\bbx^*)- \frac{1}{n}\sum_{i=1}^n F_i(\bar{\bbx}^0)\right]
+ \sum_{t=0}^{T-1}\frac{\left(1- 1/\sqrt{T}\right)^t GD}{T}
 \nonumber\\
&\qquad +\sum_{t=0}^{T-1}\frac{LD^2}{ T^{3/2}}+\sum_{t=0}^{T-1} \frac{ LD^2}{(1-\beta)T^2} 
+\sum_{t=0}^{T-1} \frac{ GD}{(1-\beta)T^{3/2}}
+\sum_{t=0}^{T-1}\frac{LD^2}{2T^2}.
\end{align}
%%%%
By using the inequality $\sum_{t=0}^{T-1}(1- 1/\sqrt{T})^t\leq \sqrt{T}$ and simplifying the terms on the right hand side \eqref{final_thm_proof_1200} we obtain that
to the expression
%%%%
\begin{align} \label{final_thm_proof_1300}
& \frac{1}{n}\sum_{i=1}^n F_i(\bbx^*)-\frac{1}{n}\sum_{i=1}^n F_i(\bar{\bbx}^{T}) \nonumber\\
&\leq\frac{1}{e} \left[ \frac{1}{n}\sum_{i=1}^n F_i(\bbx^*)- \frac{1}{n}\sum_{i=1}^n F_i(\bar{\bbx}^0)\right]
+\frac{ GD}{T^{1/2}}
+\frac{LD^2}{ T^{1/2}}  + \frac{ LD^2}{(1-\beta)T} 
+\frac{ GD}{(1-\beta)T^{1/2}}
+\frac{LD^2}{2T}\nonumber\\
&=\frac{1}{e} \left[ \frac{1}{n}\sum_{i=1}^n F_i(\bbx^*)- \frac{1}{n}\sum_{i=1}^n F_i(\bar{\bbx}^0)\right]
+\frac{ LD^2+GD(1+(1-\beta)^{-1})}{T^{1/2}}
+ \frac{ LD^2(0.5+(1-\beta)^{-1})}{T},
\end{align}
%%%%
where to derive the first inequality we used $(1-1/T)^T\leq 1/e$. Note that we set $\bbx_i^0=\bb0_p$ for all $i\in\ccalN$ and therefore $\bar{\bbx}^0=\bb0_p$. Since we assume that $F_i(\bb0_p)\geq 0$ for all $i\in\ccalN$, it implies that $\frac{1}{n}\sum_{i=1}^n F_i(\bar{\bbx}^0)=\frac{1}{n}\sum_{i=1}^n F_i(\bb0_p)\geq 0$ and the
expression in \eqref{final_thm_proof_1300} can be simplified to
%%%%
\begin{align} \label{final_thm_proof_1400}
\frac{1}{n}\sum_{i=1}^n F_i(\bar{\bbx}^T)\geq (1-e^{-1} ) \frac{1}{n}\sum_{i=1}^n F_i(\bbx^*)
-  \frac{ LD^2+GD(1+(1-\beta)^{-1})}{T^{1/2}}
- \frac{ LD^2(0.5+(1-\beta)^{-1})}{T} .
\end{align}
%%%
Also, since the norm of local gradients is uniformly bounded by $G$, the local functions $F_i$ are $G$-Lipschitz. This observation implies that
\begin{align} \label{final_thm_proof_1500}
\left|\frac{1}{n}\sum_{i=1}^n F_i(\bar{\bbx}^T)-\frac{1}{n}\sum_{i=1}^n F_i(\bbx_j^T)\right| \leq \frac{G}{n}\sum_{i=1}^n\|\bar{\bbx}^T-\bbx_j^T\| \leq \frac{GD}{T(1-\beta)},
\end{align}
%%%
where the second inequality holds by using the result in Lemma \ref{lemma:eq:bound_on_dif_from_avg} and the Cauchy-Schwartz inequality. Therefore, by combining the results in \eqref{final_thm_proof_1400} and \eqref{final_thm_proof_1500} we obtain that for all $j=\ccalN$
\begin{align} 
\frac{1}{n}\sum_{i=1}^n F_i(\bbx_j^T)
&\geq (1-e^{-1} ) \frac{1}{n}\sum_{i=1}^n F_i(\bbx^*) 
-  \frac{ LD^2+GD(1+(1-\beta)^{-1})}{T^{1/2}}\nonumber\\
&\qquad - \frac{ GD(1-\beta)^{-1}+LD^2(0.5+(1-\beta)^{-1})}{T},
\end{align}
%%%
and the claim in \eqref{local_node_bound} follows.

%%%%%%%%%%%%%%%%%%%%%%%%%%%%%%%%%%%
%%%%%%%%%%%%%%%%%%%%%%%%%%%%%%%%%%%
%%%%%%%%%%%%%%%%%%%%%%%%%%%%%%%%%%%
%%%%   S  U  B  S  E  C  T  I  O  N    %%%%%%%%%%%%%%%%
%%%%%%%%%%%%%%%%%%%%%%%%%%%%%%%%%%%
%%%%%%%%%%%%%%%%%%%%%%%%%%%%%%%%%%%

\subsection{How to Construct an Unbiased Estimator of the Gradient in Multilinear Extensions}\label{unbiased}
%Recall that $\mbE_{\bbz\sim P} [\tilde{f}(S, \bbz)]$. In terms of the multilinear extensions, we obtain $F(\bbx_{con}) = \mbE_{\bbz\sim P} [\tilde{F}(\bbx_{con}, \bbz)]$, where $F$ and $ \tilde{F}$ denote the multilinear extension for $f$ and $\tilde{f}$, respectively. So $\nabla \tilde{F}(\bbx_{con}, \bbz)$ is an unbiased estimator of $\nabla F(\bbx_{con})$ when $\bbz\sim P$. Note that $\tilde{F}(\bbx_{con}, \bbz)$ is a multilinear extension. 

In this section, we provide an unbiased estimator for the gradient of a multilinear extension. We thus consider an arbitrary submodular set function $h: 2^V \to \mathbb{R}$ with multilinear $H$. Our goal is to provide an unbiased estimator for $\nabla H(\bbx)$. We have $H(\bbx) = \sum_{S\subseteq V} \prod_{i\in S}x_i\prod_{j\not\in S} (1-x_j) h(S)$. Now, it can easily be shown that 
$$\frac{\partial H}{\partial x_i} = H(\bbx;x_i\leftarrow 1) - H(\bbx;x_i\leftarrow 0).$$
where for example by $(\bbx;x_i\leftarrow 1)$ we mean a vector which has value $1$ on its $i$-th coordinate and is equal to $\bbx$ elsewhere. To create an unbiased estimator for $\frac{\partial H}{\partial x_i} $ at a point $\bbx$ we can simply sample a set $S$ by including each element in it independently with probability $x_i$ and use $h(S \cup \{i\}) - h(S \setminus \{i\})$ as an unbiased estimator for the $i$-th partial derivative. We can sample one single set $S$ and use the above trick for all the coordinates.  This involves $n$ function computations for $h$. Having a mini-batch size $B$ we can repeat this procedure $B$ times and then average.  

Note that since every element of the unbiased estimator is of the form $h(S \cup \{i\}) - h(S \setminus \{i\})$ for some chosen set $S$, then due to submodularity of the function $h$ every element of the unbiased estimator is bounded above by the maximum marginal value of $h$ (i.e. $\max_{i \in V}$ h(\{i\})). As a result, the norm of the unbiased estimator (of the gradient of $H$) is bounded above by $\sqrt{|V|} \max_{i \in V} h(\{i\})$.

%%%%%%%%%%%%%%%%%%%%%%%%%%%%%%%%%%%
%%%%%%%%%%%%%%%%%%%%%%%%%%%%%%%%%%%
%%%%%%%%%%%%%%%%%%%%%%%%%%%%%%%%%%%
%%%%   S  U  B  S  E  C  T  I  O  N    %%%%%%%%%%%%%%%%
%%%%%%%%%%%%%%%%%%%%%%%%%%%%%%%%%%%
%%%%%%%%%%%%%%%%%%%%%%%%%%%%%%%%%%%

\subsection{Proof of Theorem \ref{theorem:main_theorem_discrete}}\label{proof:theorem:main_theorem_discrete}

The steps of the proof are similar to the one for Theorem 1. In particular, for the Discrete DCG method we can also show that the expressions in \eqref{final_thm_proof_100}-\eqref{final_thm_proof_800} hold and we can write
%%%
\begin{align} \label{final_thm_discrete_proof_100}
&\frac{1}{n}\sum_{i=1}^n F_i(\bar{\bbx}^{t+1})-\frac{1}{n}\sum_{i=1}^n F_i(\bar{\bbx}^{t})  \nonumber\\
&\geq \frac{1}{nT} \left[ \sum_{i=1}^n F_i(\bbx^*)- \sum_{i=1}^n F_i(\bar{\bbx}^t)\right]
-\frac{D}{nT}   \left\|\sum_{i=1}^n \bbd_i^t - \sum_{i=1}^n \nabla F_i(\bar{\bbx}^t) \right\|   
-\frac{D}{nT} \sum_{j=1}^n  \left\| \frac{1}{n}\sum_{i=1}^n\bbd_i^t-\bbd_j^t\right\|  
-\frac{LD^2}{2T^2}.
\end{align}
%%%%
Now we proceed to derive upper bounds for the norms on the right hand side of \eqref{final_thm_discrete_proof_100}. To derive these bounds we use the results in Lemmata~\ref{lemma:ar_in_avg_bound} and \ref{lemma:eq:bound_on_dif_from_avg} which also hold for the  Discrete DCG algorithm.

We first derive an upper bound for the sum $ \sum_{j=1}^n  \| \frac{1}{n}\sum_{i=1}^n\bbd_i^t-\bbd_j^t\|  $ in \eqref{final_thm_discrete_proof_100}. To achieve this goal the following lemma is needed. 

%%%%%%%%%%%%%%%%%%%%%%%%%%%%%%%%%%%
%%%%%%%%%%%%%%%%%%%%%%%%%%%%%%%%%%%
%%%%%%   L  E  M  M  A  %%%%%%%%%%%%%%%%%%%%%
%%%%%%%%%%%%%%%%%%%%%%%%%%%%%%%%%%%
%%%%%%%%%%%%%%%%%%%%%%%%%%%%%%%%%%%
\begin{lemma}\label{lemma:bound_on_g}
Consider the proposed Discrete DCG method defined in Algorithm \ref{algo_DDCG}. If Assumptions \ref{ass:smoothness2} and \ref{ass:bounded_variance} hold, then for all $i\in \ccalN$ and $t\geq0$ the expected squared norm $\E{\|  \bbg_i^t\|^2}$ is bounded above by
\begin{equation}\label{bound_on_g}
\E{\|  \bbg_i^t\|^2} \leq K^2,
\end{equation} 
where $K^2=\sigma^2+G^2$.
\end{lemma}

%%%%%%%%%%%%%%%%%%%%%%%%%%%%%%%%%%%
%%%%%%%%%%%%%%%%%%%%%%%%%%%%%%%%%%%
%%%%%%.  P  R  O  O  F  %%%%%%%%%%%%%%%%%%%%%
%%%%%%%%%%%%%%%%%%%%%%%%%%%%%%%%%%%
%%%%%%%%%%%%%%%%%%%%%%%%%%%%%%%%%%%
\begin{myproof}
Considering the condition in Assumption \ref{ass:bounded_variance} on the variance of stochastic gradients, we can define $K^2:=\sigma^2+G^2$ as an upper bound on the expected norm of stochastic gradients, i.e., for all $\bbx\in\ccalC$ and $i\in \ccalN$
%%%%
\begin{equation} \label{final_thm_discrete_proof_300}
\E{\|\nabla \tilde{F}_i(\bbx_i^t)\|^2} \leq K^2.
\end{equation}
%%%
Now we use an induction argument to show that the expected norm $\E{\|\bbg_{i}^t\|^2}\leq K^2$. %where $\ccalF^t$ is the sigma-algebra that measures the history of the system up until time $t$. 
Since the iterates are initialized at $\bbg_{i}^0=\bb0$, the update in \eqref{gradient_approx_update} implies that $\E{\|\bbg_{i}^1\|^2\mid \bbx_i^1}= \phi^2 \E{\|\nabla \tilde{F}_i(\bbx_i^1)\|^2 \mid \bbx_i^1}\leq \phi^2 K^2\leq K^2$. Since $\E{\E{\|\bbg_{i}^1\|^2\mid \bbx_i^1}}=\E{\|\bbg_{i}^1\|^2}$ it follows that $\E{\|\bbg_{i}^1\|^2}\leq K^2$. Now we proceed to show that if  $\E{\|\bbg_{i}^{t-1}\|^2}\leq K^2$ then  $\E{\|\bbg_{i}^t\|^2}\leq K^2$.

Recall the update of $\bbg_i^t$ in \eqref{gradient_approx_update}. By computing the squared norm of both sides and using the Cauchy-Schwartz inequality we obtain that
%%%
\begin{align}
\|\bbg_{i}^t\|^2 \leq  (1-\phi)^2 \|\bbg_{i}^{t-1}\|^2 + \phi^2 \|  \nabla \tilde{F}_i(\bbx_i^t)\|^2 + 2\phi(1-\phi)\|\bbg_{i}^{t-1}\| \|  \nabla \tilde{F}_i(\bbx_i^t)\|.
\end{align}
%%%
Compute the expectation with respect to the random variable corresponding to the stochastic gradient $\nabla \tilde{F}_i(\bbx_i^t)$ to obtain 
%%%
\begin{align}\label{sdaasdad}
\E{\|\bbg_{i}^t\|^2\mid \bbx_i^t} \leq  (1-\phi)^2 \|\bbg_{i}^{t-1}\|^2 + \phi^2 \E{\|  \nabla \tilde{F}_i(\bbx_i^t)\|^2 \mid \bbx_i^t}+ 2\phi(1-\phi)\|\bbg_{i}^{t-1}\|\E{ \|  \nabla \tilde{F}_i(\bbx_i^t)\|\mid \bbx_i^t}.
\end{align}
%%%
Note that according to Jensen's inequality $\E{\|\nabla \tilde{F}_i(\bbx_i^t)\|^2} \leq K^2$ implies that $\E{\|\nabla \tilde{F}_i(\bbx_i^t)\|} \leq K$. Replacing these bounds into \eqref{sdaasdad} yields 
%%%
\begin{align}
\E{\|\bbg_{i}^t\|^2\mid \bbx_i^t} \leq  (1-\phi)^2 \|\bbg_{i}^{t-1}\|^2 + \phi^2K^2+ 2K\phi(1-\phi)\|\bbg_{i}^{t-1}\|.
\end{align}
%%%
Now by computing the expectation of both sides with respect to all sources of randomness from $t=0$ and using the simplification $\E{\E{\|\bbg_{i}^t\|^2\mid \bbx_i^t}}=\E{\|\bbg_{i}^t\|^2}$ we can write 
%%%
\begin{align}
\E{\|\bbg_{i}^t\|^2\mid }
& \leq  (1-\phi)^2 \E{\|\bbg_{i}^{t-1}\|^2} + \phi^2K^2+ 2K\phi(1-\phi)\E{\|\bbg_{i}^{t-1}\|}\nonumber\\
&\leq  (1-\phi)^2 K^2+ \phi^2K^2+ 2K\phi(1-\phi)K 
\nonumber\\
&= K^2,
\end{align}
%%%
and the claim in \eqref{bound_on_g} follows by induction. 
\end{myproof}
%
%\red{
%Based on the update in \eqref{gradient_approx_update}, if we set $\bbg_i^0=0$, then for $t\geq 1$ the iterate $\bbg_i^t$ can be written as
%%%%
%\begin{align}
%\bbg_i^t = \sum_{s=1}^t \phi(1-\phi)^{t-s}  \nabla \tilde{F}_i(\bbx_i^s).
%\end{align}
%%%%
%Evaluate norm of both sides and used the Cauchy-Schwartz inequality to obtain
%%%%
%\begin{align}
%\|\bbg_i^t \|\leq \sum_{s=1}^t \phi(1-\phi)^{t-s} \left\| \nabla \tilde{F}_i(\bbx_i^s) \right\|.
%\end{align}
%%%%
%Computing the expected value of both sides and using the inequality $\E{\| \nabla \tilde{F}_i(\bbx_i^s) \|}\leq \left(\E{\| \nabla \tilde{F}_i(\bbx_i^s) \|^2}\right)^{1/2} \leq K$ lead to
%%%%
%\begin{align}\label{plplpl}
%\E{\|\bbg_i^t \|} &\leq \sum_{s=1}^t \phi(1-\phi)^{t-s}K.
%\end{align}
%%%%
%}
%The claim in \eqref{bound_on_g} follows by simplifying the sum on the right hand side of \eqref{plplpl}.

%%%%%%%%%%%%%%%%%%%%%%%%%%%%%%%%%%%
%%%%%%%%%%%%%%%%%%%%%%%%%%%%%%%%%%%
%%%%%%. MAIN MATTER  %%%%%%%%%%%%%%%%%%%%
%%%%%%%%%%%%%%%%%%%%%%%%%%%%%%%%%%%
%%%%%%%%%%%%%%%%%%%%%%%%%%%%%%%%%%%
We use the result in Lemma \ref{lemma:bound_on_g} to find an upper bound for the sum $(1/n)\sum_{j=1}^n  \left\| \bar{\bbd}^t-\bbd_j^t\right\|  $ on the right hand side of \eqref{final_thm_discrete_proof_100}.

%%%%%%%%%%%%%%%%%%%%%%%%%%%%%%%%%%%
%%%%%%%%%%%%%%%%%%%%%%%%%%%%%%%%%%%
%%%%%%.  L  E  M  M  A  %%%%%%%%%%%%%%%%%%%%%
%%%%%%%%%%%%%%%%%%%%%%%%%%%%%%%%%%%
%%%%%%%%%%%%%%%%%%%%%%%%%%%%%%%%%%%
\begin{lemma}\label{lemma:bound_on_d}
Consider the proposed Discrete DCG method defined in Algorithm \ref{algo_DDCG}. If Assumptions \ref{ass:weights}, \ref{ass:smoothness2} and \ref{ass:bounded_variance} hold, then for all $i\in \ccalN$ and $t\geq0$ we have
\begin{equation}\label{bound_on_d}
\E{\frac{1}{n}\sum_{i=1}^n\|\bbd_i^t-\bar{\bbd}^t\|}\leq \frac{\alpha K}{1-\beta(1-\alpha)},
\end{equation} 
where $K=(\sigma^2+G^2)^{1/2}$.
\end{lemma}

%%%%%%%%%%%%%%%%%%%%%%%%%%%%%%%%%%%
%%%%%%%%%%%%%%%%%%%%%%%%%%%%%%%%%%%
%%%%%%.  P  R  O  O  F  %%%%%%%%%%%%%%%%%%%%%
%%%%%%%%%%%%%%%%%%%%%%%%%%%%%%%%%%%
%%%%%%%%%%%%%%%%%%%%%%%%%%%%%%%%%%%
\begin{myproof}
Define the vector $\bbg_{con}^t=[\bbg_1^t;\dots;\bbg_n^t]$ as the concatenation of the local vectors $\bbg_i^t$ at time $t$. Further, recall the definitions of the vectors $\bbx_{con}=[\bbx_1;\dots;\bbx_n]\in \reals^{np}$ and $\bbd_{con}=[\bbd_1;\dots;\bbd_n]\in \reals^{np}$ as the concatenation of the local variables and local approximate gradients, respectively, and the definition of $\bar{\bbd}_{con}^t=[\bar{\bbd}^t;\dots;\bar{\bbd}^t]$ as the concatenation of $n$ copies of the average vector $\bar{\bbd}^t$. By following the steps of the proof for Lemma~\ref{lemma:bound_on_gradient_consensus_error}, it can be shown that 
%%%%
\begin{align}\label{final_thm_discrete_proof_500}
\|\bbd_{con}^t -\bar{\bbd}_{con}^t\| &=  \left\| \alpha \sum_{s=1}^{t}(1-\alpha)^{t-s} \left[(\bbW^{t-s}-\frac{\bbone_n\bbone_n^{\dag}}{n})\otimes \bbI \right] \bbg_{con}^t  \right\|  \nonumber\\
&\leq   \alpha \sum_{s=1}^{t}(1-\alpha)^{t-s} \left\|(\bbW^{t-s}-\frac{\bbone_n\bbone_n^{\dag}}{n})\otimes \bbI \right\| \left\| \bbg_{con}^t\right\|  \nonumber\\
&\leq   \alpha \sum_{s=1}^{t}(1-\alpha)^{t-s} \beta^{t-s}\left\| \bbg_{con}^t\right\|  .
\end{align}
%%%%
By computing the expected value of both sides and using the result in \eqref{bound_on_g} we obtain that
%%%%
\begin{align}\label{final_thm_discrete_proof_500}
\E{\|\bbd_{con}^t -\bar{\bbd}_{con}^t\|} 
&\leq   \alpha \sqrt{n}K \sum_{s=1}^{t}(1-\alpha)^{t-s} \beta^{t-s}  \nonumber\\
&\leq   \frac{\alpha \sqrt{n} K}{1-\beta(1-\alpha)},
\end{align}
%%%%
where to derive the first inequality we use the fact that 
\begin{equation}
\E{\| \bbg_{con}^t\| }\leq(\E{\| \bbg_{con}^t\|^2})^{1/2}=(\E{(\sum_{i=1}^n \|\bbg_i^t\|^2)})^{1/2}=(\sum_{i=1}^n\E{ \|\bbg_i^t\|^2})^{1/2}\leq \sqrt{n}K.
\end{equation}
By combining the result in \eqref{final_thm_discrete_proof_500} with the inequality 
\begin{equation}
\frac{1}{n}\sum_{i=1}^n\|\bbd_i^t-\bar{\bbd}^t\|\leq \frac{1}{\sqrt{n}} \left[\sum_{i=1}^n\|\bbd_i^t-\bar{\bbd}^t\|^2\right]^{1/2} =\frac{1}{\sqrt{n}}  \|\bbd_{con}^t -\bar{\bbd}_{con}^t\|,
\end{equation}
the claim in \eqref{bound_on_d} follows.
\end{myproof}
%%%%

%%%%%%%%%%%%%%%%%%%%%%%%%%%%%%%%%%%
%%%%%%%%%%%%%%%%%%%%%%%%%%%%%%%%%%%
%%%%%%. MAIN MATTER  %%%%%%%%%%%%%%%%%%%%
%%%%%%%%%%%%%%%%%%%%%%%%%%%%%%%%%%%
%%%%%%%%%%%%%%%%%%%%%%%%%%%%%%%%%%%
The result in Lemma \ref{lemma:bound_on_d} shows $\frac{1}{n}\sum_{i=1}^n\|\bbd_i^t-\bar{\bbd}^t\| $ is bounded above by $(\alpha K)/({1-\beta(1-\alpha)})$ in expectation.  To bound the second sum in \eqref{final_thm_discrete_proof_100}, which is$ \left\|\sum_{i=1}^n \bbd_i^t - \sum_{i=1}^n \nabla F_i(\bar{\bbx}^t) \right\|  $, we first introduce the following lemma, which was presented in \citep{mokhtari2017conditional} in a slightly different form.

%%%%%%%%%%%%%%%%%%%%%%%%%%%%%%%%%%%
%%%%%%%%%%%%%%%%%%%%%%%%%%%%%%%%%%%
%%%%%%.  L  E  M  M  A  %%%%%%%%%%%%%%%%%%%%%
%%%%%%%%%%%%%%%%%%%%%%%%%%%%%%%%%%%
%%%%%%%%%%%%%%%%%%%%%%%%%%%%%%%%%%%
\begin{lemma}\label{lemma:bound_on_stoc}
Consider the proposed Discrete DCG method defined in Algorithm \ref{algo_DDCG}. If Assumptions \ref{ass:weights}-\ref{ass:bounded_variance} hold and we set $\phi=1/T^{2/3}$, then for all $i\in \ccalN$ and $t\geq0$ we have
\begin{align}\label{bound_on_stoch}
\E{\sum_{i=1}^n\|\nabla F_i(\bbx_i^t) - \bbg_i^t\|^2}
      &\leq 
   \left(1-\frac{1}{2T^{2/3}}\right)^tnG^2 
   + \frac{6nL^2D^2C}{T^{4/3}}
   +\frac{2n\sigma^2+12nL^2D^2C}{T^{2/3}},
\end{align}
where $C:=1+({2}/{(1-\beta)^2})$.
\end{lemma}

%%%%%%%%%%%%%%%%%%%%%%%%%%%%%%%%%%%
%%%%%%%%%%%%%%%%%%%%%%%%%%%%%%%%%%%
%%%%%%.  P  R  O  O  F  %%%%%%%%%%%%%%%%%%%%%
%%%%%%%%%%%%%%%%%%%%%%%%%%%%%%%%%%%
%%%%%%%%%%%%%%%%%%%%%%%%%%%%%%%%%%%
\begin{myproof}
Use the update $\bbg_i^t := (1-\phi) \bbg_i^{t-1} + \phi \nabla \tilde{F}(\bbx_i^t)$  to write  the squared norm $\|\nabla F_i(\bbx_i^t) - \bbg_i^t\|^2$ as
%%%%
\begin{align}\label{proof:bound_on_grad_100}
\|\nabla F_i(\bbx_i^t) -  \bbg_i^t\|^2 
=\|\nabla F_i(\bbx_i^t) -(1-\phi) \bbd_{t-1} - \phi \nabla \tilde{F}_i(\bbx_i^t)\|^2.
\end{align}
%%%
Add and subtract the term $(1-\phi)\nabla F_i(\bbx_i^{t-1})$ to the right hand side of \eqref{proof:bound_on_grad_100} and regroup the terms to obtain
%%%%
\begin{align}\label{proof:bound_on_grad_200}
&\|\nabla F_i(\bbx_i^t) - \bbg_i^t\|^2
\nonumber\\
&=\|\phi(\nabla F_i(\bbx_i^t)-\nabla \tilde{F}_i(\bbx_i^t))  
+(1-\phi)(\nabla F_i(\bbx_i^t)-\nabla F_i(\bbx_i^{t-1}))
+(1-\phi)(\nabla F_i(\bbx_i^{t-1}) - \bbg_i^{t-1} )\|^2.
\end{align}
%%%%
Define $\ccalF^t$ as a sigma algebra that measures the history of the system up until time $t$. Expanding the square and computing the conditional expectation $\E{\cdot\mid \ccalF^t}$ of the resulted expression yield 
%%%
\begin{align}\label{proof:bound_on_grad_300}
&\E{\|\nabla F_i(\bbx_i^t) - \bbg_i^t\|^2\mid\ccalF^t} =\phi^2\E{\|\nabla F_i(\bbx_i^t)-\nabla \tilde{F}_i(\bbx_i^t)\|^2\mid\ccalF^t}+   (1-\phi)^2\|\nabla F_i(\bbx_i^{t-1}) - \bbg_i^{t-1} \|^2 
     \nonumber\\
     &  \quad
   +(1-\phi)^2\|\nabla F_i(\bbx_i^t)-\nabla F_i(\bbx_i^{t-1})\|^2
    + 2 (1-\phi)^2\langle \nabla F_i(\bbx_i^t)-\nabla F_i(\bbx_i^{t-1}) , \nabla F_i(\bbx_i^{t-1}) -  \bbg_i^{t-1} \rangle,
 \end{align}
%%%%
where we have used the fact $\E{\nabla \tilde{F}_i(\bbx_i^t) \mid\ccalF^t}=\nabla F_i(\bbx_i^t) $. 
  The term $\E{\|\nabla F_i(\bbx_i^t)-\nabla \tilde{F}_i(\bbx_i^t)\|^2\mid\ccalF^t}$ can be bounded above by $\sigma^2$ according to Assumption \ref{ass:bounded_variance}. Based on Assumption~\ref{ass:smoothness}, we can also show that the squared norm $\|\nabla F_i(\bbx_i^t)-\nabla F_i(\bbx_i^{t-1})\|^2$ is upper bounded by $L^2\|\bbx_i^t-\bbx_i^{t-1}\|^2$. Moreover, the inner product $2\langle \nabla F_i(\bbx_i^t)\!-\!\nabla F_i(\bbx_i^{t-1}) , \nabla F_i(\bbx_i^{t-1}) - \bbd_{t-1} \rangle$ can be upper bounded by $\zeta \|\nabla F_i(\bbx_i^{t-1}) - \bbd_{t-1}\|^2+(1/\zeta) L^2\|\bbx_i^t-\bbx_i^{t-1}\|^2$ using Young's inequality (i.e., $2\langle \bba,\bbb\rangle \leq \zeta\|\bba\|^2+\|\bbb\|^2/\beta$  for any $\bba,\bbb\in \reals^n$ and $\zeta>0$) and the condition in Assumption \ref{ass:smoothness}, where $\zeta>0$ is a free scalar. Applying these substitutions into \eqref{proof:bound_on_grad_300} leads to 
  %%%
\begin{align}\label{proof:bound_on_grad_400}
&\E{\|\nabla F_i(\bbx_i^t) -  \bbg_i^t\|^2\mid\ccalF^t}
\nonumber\\
& \leq \phi^2\sigma^2
   +(1-\phi)^2 (1+\zeta^{-1})L^2\|\bbx_i^t-\bbx_i^{t-1}\|^2
   +(1-\phi)^2(1+\zeta)\|\nabla F_i(\bbx_i^{t-1}) -  \bbg_i^{t-1} \|^2. 
   \end{align}
By setting $\zeta=\phi/2$ we can replace $(1-\phi)^2 (1+\zeta^{-1})$ and $(1-\phi)^2(1+\zeta)$ by their upper bounds $(1+2\phi^{-1})$ and $(1-\phi/2)$, respectively. Applying theses substitutions and summing up both sides of the resulted inequality for $i=1,\dots,n$ lead to
  %%%
\begin{align}\label{proof:bound_on_grad_401}
&\E{\sum_{i=1}^n\|\nabla F_i(\bbx_i^t) - \bbg_i^t\|^2\mid\ccalF^t}
\nonumber\\
&
\leq n\phi^2\sigma^2
   + L^2(1+2\phi^{-1}) \sum_{i=1}^n\|\bbx_i^t-\bbx_i^{t-1}\|^2
   +\left(1-\frac{\phi}{2}\right)\sum_{i=1}^n\|\nabla F_i(\bbx_i^{t-1}) - \bbg_i^{t-1}\|^2.
   \end{align}
Now we proceed to derive an upper bound for the sum $\sum_{i=1}^n\|\bbx_i^t-\bbx_i^{t-1}\|^2$. Note that using the Cauchy-Schwartz inequality and the results in Lemmata~\ref{lemma:ar_in_avg_bound} and \ref{lemma:eq:bound_on_dif_from_avg} we can show that
%%%
\begin{align}\label{proof:bound_on_grad_402}
 \sum_{i=1}^n\|\bbx_i^t-\bbx_i^{t-1}\|^2
&\leq \sum_{i=1}^n\left(3 \left\|\bbx_i^t-\bar{\bbx}^t\right\|^2 +3 \left\|\bar{\bbx}^t-\bar{\bbx}^{t-1}\right\|^2+3\left\|\bar{\bbx}^{t-1}-{\bbx}_i^{t-1}\right\|^2\right)\nonumber\\
&\leq \frac{3nD^2}{T^2(1-\beta)^2} + \frac{3nD^2}{T^2}+\frac{3nD^2}{T^2(1-\beta)^2}\nonumber\\
&=  \frac{3nD^2}{T^2}\left(1+\frac{2}{(1-\beta)^2}\right).
\end{align}
%%%
Replace the sum $ \sum_{i=1}^n\|\bbx_i^t-\bbx_i^{t-1}\|^2$ in \eqref{proof:bound_on_grad_401} by its upper bound in \eqref{proof:bound_on_grad_402} and compute the expectation with respect to $\ccalF_0$ to obtain
%%%
\begin{align}\label{proof:bound_on_grad_800}
&\E{\sum_{i=1}^n\|\nabla F_i(\bbx_i^t) - \bbg_i^t\|^2}\nonumber\\
&\leq 
   \left(1-\frac{\phi}{2}\right)\E{\sum_{i=1}^n\|\nabla F_i(\bbx_i^{t-1}) - \bbg_i^{t-1} \|^2}
   +n\phi^2\sigma^2
   + (1+2\phi^{-1})\frac{3nL^2D^2}{T^2}\left(1+\frac{2}{(1-\beta)^2}\right).
\end{align}
%%%
Set $\phi=T^{-2/3}$ to obtain 
%%%
\begin{align}\label{proof:bound_on_grad_900}
&\E{\sum_{i=1}^n\|\nabla F_i(\bbx_i^t) - \bbg_i^t\|^2}\nonumber\\
&\leq 
   \left(1-\frac{1}{2T^{2/3}}\right)\E{\sum_{i=1}^n\|\nabla F_i(\bbx_i^{t-1}) - \bbg_i^{t-1} \|^2}
   +\frac{n\sigma^2}{T^{4/3}}
   + \frac{3nL^2D^2C}{T^2}
   +\frac{6nL^2D^2C}{T^{4/3}},
\end{align}
%%%
where $C:=\left(1+\frac{2}{(1-\beta)^2}\right)$. Applying the expression in \eqref{proof:bound_on_grad_900} recursively leads to 
%%%
\begin{align}\label{proof:bound_on_grad_1000}
&\E{\sum_{i=1}^n\|\nabla F_i(\bbx_i^t) - \bbg_i^t\|^2}\nonumber\\
&\leq 
   \left(1-\frac{1}{2T^{2/3}}\right)^t\sum_{i=1}^n\|\nabla F_i(\bbx_i^{0}) - \bbd_{0} \|^2  +\left(\frac{n\sigma^2}{T^{4/3}}
   + \frac{3nL^2D^2C}{T^2}
   +\frac{6nL^2D^2C}{T^{4/3}}\right) \sum_{s=0}^{t-1} \left(1-\frac{1}{2T^{2/3}}\right)^s\nonumber\\
   &\leq 
   \left(1-\frac{1}{2T^{2/3}}\right)^t\sum_{i=1}^n\|\nabla F_i(\bbx_i^{0}) - \bbd_{0} \|^2  +\frac{2n\sigma^2}{T^{2/3}}
   + \frac{6nL^2D^2C}{T^{4/3}}
   +\frac{12nL^2D^2C}{T^{2/3}}\nonumber\\
      &\leq 
   \left(1-\frac{1}{2T^{2/3}}\right)^tnG^2  +\frac{2n\sigma^2}{T^{2/3}}
   + \frac{6nL^2D^2C}{T^{4/3}}
   +\frac{12nL^2D^2C}{T^{2/3}},
\end{align}
%%%
and the claim in \eqref{bound_on_stoch} follows. 
\end{myproof}

%%%%%%%%%%%%%%%%%%%%%%%%%%%%%%%%%%%
%%%%%%%%%%%%%%%%%%%%%%%%%%%%%%%%%%%
%%%%%%. MAIN MATTER  %%%%%%%%%%%%%%%%%%%%
%%%%%%%%%%%%%%%%%%%%%%%%%%%%%%%%%%%
%%%%%%%%%%%%%%%%%%%%%%%%%%%%%%%%%%%
We use the result in Lemma \ref{lemma:bound_on_stoc} to derive an upper bound for $\left\|\frac{1}{n}\sum_{i=1}^n \bbd_i^t  - \frac{1}{n}\sum_{i=1}^n\nabla F_i(\bar{\bbx}^t) \right\|$ in expectation. 

%%%%%%%%%%%%%%%%%%%%%%%%%%%%%%%%%%%
%%%%%%%%%%%%%%%%%%%%%%%%%%%%%%%%%%%
%%%%%%.  L  E  M  M  A  %%%%%%%%%%%%%%%%%%%%%
%%%%%%%%%%%%%%%%%%%%%%%%%%%%%%%%%%%
%%%%%%%%%%%%%%%%%%%%%%%%%%%%%%%%%%%
\begin{lemma}\label{lemma:bound_on_whatever}
Consider the proposed Discrete DCG method defined in Algorithm \ref{algo_DDCG}. If Assumptions \ref{ass:weights}-\ref{ass:bounded_variance} hold and we set $\alpha= 1/\sqrt{T}$ and $\phi=1/T^{2/3}$, then for all $i\in \ccalN$ and $t\geq0$ we have
\begin{align}\label{bound_on_whatever}
\E{\left\|\frac{1}{n}\sum_{i=1}^n \bbd_i^t  - \frac{1}{n}\sum_{i=1}^n\nabla F_i(\bar{\bbx}^t) \right\|}
&  \leq G\left(1-\frac{1}{T^{1/2}}\right)^t 
 +G\left(1-\frac{1}{2T^{2/3}}\right)^{t/2}  
 +\frac{LD}{ T^{1/2}} \nonumber\\
&\qquad + \frac{ LD}{T(1-\beta)}
   + \frac{\sqrt{6} LDC^{1/2}}{T^{2/3}}
   +\frac{\sqrt{2}\sigma+\sqrt{12} LDC^{1/2}}{T^{1/3}},
\end{align}
where $C:=1+({2}/{(1-\beta)^2})$.
\end{lemma}

%%%%%%%%%%%%%%%%%%%%%%%%%%%%%%%%%%%
%%%%%%%%%%%%%%%%%%%%%%%%%%%%%%%%%%%
%%%%%%.  P  R  O  O  F  %%%%%%%%%%%%%%%%%%%%%
%%%%%%%%%%%%%%%%%%%%%%%%%%%%%%%%%%%
%%%%%%%%%%%%%%%%%%%%%%%%%%%%%%%%%%%
\begin{myproof}
The steps of this proof are similar to the ones in the proof of Lemma \ref{lemma:bound_on_gradient_error}. It can be shown that
%%%
\begin{align}\label{a}
&\left\|\sum_{i=1}^n \bbd_i^t  -\sum_{i=1}^n\nabla F_i(\bar{\bbx}^t)\right\|\nonumber\\
&=\left\| (1-\alpha)  \sum_{j=1}^n \bbd_{j}^{t-1}  +\alpha \sum_{i=1}^n \bbg_i^t  - \sum_{i=1}^n\nabla F_i(\bar{\bbx}^t) \right\| \nonumber\\
&=\left\| (1-\alpha)  \sum_{j=1}^n \bbd_{j}^{t-1}  -(1-\alpha ) \sum_{i=1}^n \nabla F_i(\bar{\bbx}^{t-1})  + (1-\alpha)\sum_{i=1}^n\nabla F_i(\bar{\bbx}^{t-1})+\alpha \sum_{i=1}^n \bbg_i^t  - \sum_{i=1}^n\nabla F_i(\bar{\bbx}^t)  \right\| \nonumber\\
&=\left\| (1-\alpha) \left[ \sum_{j=1}^n \bbd_{j}^{t-1}-\sum_{i=1}^n\nabla F_i(\bar{\bbx}^{t-1})  \right]
 +(1-\alpha ) \left[\sum_{i=1}^n \nabla F_i(\bar{\bbx}^{t-1})  - \sum_{i=1}^n \nabla F_i(\bar{\bbx}^{t}) \right]
 +\alpha \left[\sum_{i=1}^n   \bbg_i^t   - \sum_{i=1}^n\nabla F_i(\bar{\bbx}^{t})   \right]
 \right\| \nonumber\\
 &\leq (1-\alpha) \left\| \sum_{j=1}^n \bbd_{j}^{t-1}-\sum_{i=1}^n\nabla F_i(\bar{\bbx}^{t-1})  \right\|
 +(1-\alpha ) \left\|\sum_{i=1}^n \nabla F_i(\bar{\bbx}^{t-1})  - \sum_{i=1}^n \nabla F_i(\bar{\bbx}^{t}) \right\|
 +\alpha \left\|\sum_{i=1}^n  \bbg_i^t   - \sum_{i=1}^n\nabla F_i(\bar{\bbx}^{t})  \right\| \nonumber\\
  &\leq (1-\alpha) \left\| \sum_{j=1}^n \bbd_{j}^{t-1}-\sum_{i=1}^n\nabla F_i(\bar{\bbx}^{t-1})  \right\|
 +(1-\alpha ) \left\|\sum_{i=1}^n \nabla F_i(\bar{\bbx}^{t-1})  - \sum_{i=1}^n \nabla F_i(\bar{\bbx}^{t}) \right\|
 +\alpha \left\|\sum_{i=1}^n  \bbg_i^t   - \sum_{i=1}^n \nabla F_{i}(\bbx_i^t) \right\| \nonumber\\
 &\qquad +\alpha \left\|\sum_{i=1}^n \nabla F_{i}(\bbx_i^t)  - \sum_{i=1}^n\nabla F_i(\bar{\bbx}^{t})  \right\|.
\end{align}
The first equality is the outcome of replacing $\sum_{i=1}^n \bbd_i^t $ by the expression in \eqref{app_proof_100}, the second equality is obtained by adding and subtracting $(1-\alpha ) \sum_{i=1}^n \nabla F_i(\bar{\bbx}^{t-1})$, in the third equality we regroup the terms, and the inequality follows from applying the triangle inequality twice. Applying the Cauchy--Schwarz inequality to the second and third summands in \eqref{app_proof_200} and using the Lipschitz continuity of the gradients lead to 
%%%
\begin{align}\label{app}
\left\|\sum_{i=1}^n \bbd_i^t  - \sum_{i=1}^n\nabla F_i(\bar{\bbx}^t) \right\|
&\leq (1-\alpha) \left\|  \sum_{j=1}^n \bbd_{j}^{t-1}-\sum_{i=1}^n\nabla F_i(\bar{\bbx}^{t-1}) \right\| 
+(1-\alpha ) L \sum_{i=1}^n \left\|\bar{\bbx}^{t-1} - \bar{\bbx}^{t}\right\| 
\nonumber\\
&\qquad + \alpha L \sum_{i=1}^n \|  \bbx_i^t  - \bar{\bbx}^{t}\| +\alpha \left\|\sum_{i=1}^n  \bbg_i^t   -\sum_{i=1}^n \nabla F_{i}(\bbx_i^t) \right\|\nonumber\\
&\leq (1-\alpha) \left\|  \sum_{j=1}^n \bbd_{j}^{t-1}-\sum_{i=1}^n\nabla F_i(\bar{\bbx}^{t-1}) \right\| 
+\frac{(1-\alpha )LnD}{T} \nonumber\\
&\qquad 
+ \frac{\alpha LnD}{T(1-\beta)}  +\alpha\sum_{i=1}^n \left\| \bbg_i^t   -  \nabla F_{i}(\bbx_i^t) \right\|,
\end{align}
%%%
where the last inequality follows from Lemmata~\ref{lemma:ar_in_avg_bound} and \ref{lemma:eq:bound_on_dif_from_avg}. Using the inequality 
\begin{align}
\frac{1}{\sqrt{n}} \E{\sum_{i=1}^n \left\| \bbg_i^t   -  \nabla F_{i}(\bbx_i^t) \right\|}
\leq \E{\left(\sum_{i=1}^n \left\| \bbg_i^t   -  \nabla F_{i}(\bbx_i^t) \right\|^2\right)^{1/2}}
\leq \left(\E{\sum_{i=1}^n \left\| \bbg_i^t   -  \nabla F_{i}(\bbx_i^t) \right\|^2}\right)^{1/2},
\end{align}
%%%
and the result in Lemma \ref{lemma:bound_on_whatever} we obtain that
%%%
\begin{align}\label{appppp}
 \E{\sum_{i=1}^n \left\| \bbg_i^t   -  \nabla F_{i}(\bbx_i^t) \right\|}
&\leq \sqrt{n} \left[   \left(1-\frac{1}{2T^{2/3}}\right)^tnG^2  +\frac{2n\sigma^2}{T^{2/3}}
   + \frac{6nL^2D^2C}{T^{4/3}}
   +\frac{12nL^2D^2C}{T^{2/3}}\right]^{1/2}\nonumber\\
  &\leq  nG\left(1-\frac{1}{2T^{2/3}}\right)^{t/2}   +\frac{\sqrt{2}n\sigma}{T^{1/3}}
   + \frac{\sqrt{6}nLDC^{1/2}}{T^{2/3}}
   +\frac{\sqrt{12}nLDC^{1/2}}{T^{1/3}},
\end{align}
%%%%
where the second inequality holds since $\sum_{i}a_i^2 \leq(\sum_{i}a_i)^2 $ for $a_i\geq0$. Compute the expected value of both sides of \eqref{app} and replace $ \E{\sum_{i=1}^n \left\| \bbg_i^t   -  \nabla F_{i}(\bbx_i^t) \right\|}$ by its upper bound in \eqref{appppp} to obtain
%%%
\begin{align}\label{app2222}
\E{\left\|\sum_{i=1}^n \bbd_i^t  - \sum_{i=1}^n\nabla F_i(\bar{\bbx}^t) \right\|}
&\leq (1-\alpha) \E{\left\|  \sum_{j=1}^n \bbd_{j}^{t-1}-\sum_{i=1}^n\nabla F_i(\bar{\bbx}^{t-1}) \right\| }
+\frac{(1-\alpha )LnD}{T} 
+ \frac{\alpha LnD}{T(1-\beta)}  \nonumber\\
& \quad +\alpha nG\left(1-\frac{1}{2T^{2/3}}\right)^{t/2} \!\!
   + \frac{\sqrt{6}\alpha nLDC^{1/2}}{T^{2/3}}
   +\frac{\sqrt{2}n\alpha\sigma+\sqrt{12}\alpha nLDC^{1/2}}{T^{1/3}}.
\end{align}
%%%

By multiplying both of sides of \eqref{app_proof_400} by $1/n$ and applying the resulted inequality recessively for $t$ steps we obtain 
%%%
\begin{align}\label{app_pro}
&\E{\left\|\frac{1}{n}\sum_{i=1}^n \bbd_i^t  - \frac{1}{n}\sum_{i=1}^n\nabla F_i(\bar{\bbx}^t) \right\|}\nonumber\\
&\leq (1-\alpha)^t \left\| \frac{1}{n} \sum_{j=1}^n \bbd_{j}^{0}-\frac{1}{n}\sum_{i=1}^n\nabla F_i(\bar{\bbx}^{0}) \right\| \nonumber\\
&\  + \!\left[\frac{(1\!-\!\alpha )LD}{T} 
+ \frac{\alpha LD}{T(1\!-\!\beta)} + \alpha G\left[1-\frac{1}{2T^{2/3}}\right]^{t/2}  \!\!\!
   + \frac{\sqrt{6}\alpha LDC^{1/2}}{T^{2/3}}
   +\frac{\sqrt{2}\alpha \sigma\!+\!\sqrt{12}\alpha LDC^{1/2}}{T^{1/3}}\right]\sum_{s=0}^{t-1} (1-\alpha)^s\nonumber\\
&\leq (1-\alpha)^t \frac{1}{n}\sum_{i=1}^n\left\|\nabla F_i(\bar{\bbx}^{0}) \right\| 
+\frac{(1-\alpha)LD}{\alpha T} 
+ \frac{ LD}{T(1-\beta)} +G\left(1-\frac{1}{2T^{2/3}}\right)^{t/2}  
\nonumber\\
&\qquad
   + \frac{\sqrt{6} LDC^{1/2}}{T^{2/3}}
   +\frac{\sqrt{2}\sigma+\sqrt{12} LDC^{1/2}}{T^{1/3}}\nonumber\\
&\leq \left(1-\frac{1}{T^{1/2}}\right)^t G
+\frac{(1-\alpha)LD}{\alpha T} 
+ \frac{ LD}{T(1-\beta)} +G\left(1-\frac{1}{2T^{2/3}}\right)^{t/2}  \nonumber\\
&\qquad    + \frac{\sqrt{6} LDC^{1/2}}{T^{2/3}}
   +\frac{\sqrt{2}\sigma+\sqrt{12} LDC^{1/2}}{T^{1/3}},
\end{align}
which yields the claim in \eqref{bound_on_whatever}.
\end{myproof}

%%%%%%%%%%%%%%%%%%%%%%%%%%%%%%%%%%%
%%%%%%%%%%%%%%%%%%%%%%%%%%%%%%%%%%%
%%%%%%. MAIN MATTER  %%%%%%%%%%%%%%%%%%%%
%%%%%%%%%%%%%%%%%%%%%%%%%%%%%%%%%%%
%%%%%%%%%%%%%%%%%%%%%%%%%%%%%%%%%%%
Now we can complete the proof of Theorem \ref{theorem:main_theorem_discrete} using the results in Lemmata \ref{lemma:bound_on_d} and \ref{lemma:bound_on_whatever} as well as the expression in \eqref{final_thm_discrete_proof_100}. Replace the terms on the right hand side of \eqref{final_thm_discrete_proof_100} by their upper bounds in Lemmata \ref{lemma:bound_on_d} and \ref{lemma:bound_on_whatever} to obtain
%%%
\begin{align} \label{final_thm_discrete_proof_1001}
&\E{\frac{1}{n}\sum_{i=1}^n F_i(\bar{\bbx}^{t+1})-\frac{1}{n}\sum_{i=1}^n F_i(\bar{\bbx}^{t}) } \nonumber\\
&\geq \E{ \frac{1}{nT} \left[ \sum_{i=1}^n F_i(\bbx^*)- \sum_{i=1}^n F_i(\bar{\bbx}^t)\right]}
-\left(1-\frac{1}{T^{1/2}}\right)^t \frac{DG}{T}
-\frac{(1-\alpha)LD^2}{\alpha T^2} 
-\frac{ LD^2}{T^2(1-\beta)}  \nonumber\\
&  \quad -\frac{DG}{T}\left(1-\frac{1}{2T^{2/3}}\right)^{t/2}  - \frac{\sqrt{6} LD^2C^{1/2}}{T^{5/3}}
   -\frac{\sqrt{2}\sigma+\sqrt{12} LD^2C^{1/2}}{T^{4/3}}- \frac{ D(\sigma^2+G^2)^{1/2}}{T^{3/2}(1-\beta(1-\alpha))}-\frac{LD^2}{2T^2}.
\end{align}
%%%
Regrouping the terms implies that
%%%
\begin{align} \label{final_thm_discrete_proof_10011}
&\E{ \frac{1}{n}\sum_{i=1}^n F_i(\bbx^*)-\frac{1}{n}\sum_{i=1}^n F_i(\bar{\bbx}^{t+1})  } \nonumber\\
&\leq \left(1- \frac{1}{T}\right) \E{ \frac{1}{n}\sum_{i=1}^n F_i(\bbx^*)- \frac{1}{n}\sum_{i=1}^n F_i(\bar{\bbx}^t)}
+\left(1-\frac{1}{T^{1/2}}\right)^t \frac{DG}{T}
+\frac{(1-\alpha)LD^2}{\alpha T^2} 
+\frac{ LD^2}{T^2(1-\beta)}   \nonumber\\
&  \quad  +\frac{DG}{T}\left(1-\frac{1}{2T^{2/3}}\right)^{t/2}+\frac{\sqrt{6} LD^2C^{1/2}}{T^{5/3}}
   +\frac{\sqrt{2}\sigma+\sqrt{12} LD^2C^{1/2}}{T^{4/3}}+ \frac{ D(\sigma^2+G^2)^{1/2}}{T^{3/2}(1-\beta(1-\alpha))}+\frac{LD^2}{2T^2}.
\end{align}
%%%
Now apply the expression in \eqref{final_thm_discrete_proof_10011} for $t=0,\dots,T-1$ to obtain
%%%
\begin{align} \label{final_thm_discrete_proof_10012}
&\E{ \frac{1}{n}\sum_{i=1}^n F_i(\bbx^*)-\frac{1}{n}\sum_{i=1}^n F_i(\bar{\bbx}^{T})  } \nonumber\\
&\leq \left(1- \frac{1}{T}\right)^T \E{ \frac{1}{n}\sum_{i=1}^n F_i(\bbx^*)- \frac{1}{n}\sum_{i=1}^n F_i(\bar{\bbx}^0)}
+\frac{(1-\alpha)LD^2}{\alpha T} 
+\frac{ LD^2}{T(1-\beta)}    
   +\frac{\sqrt{2}\sigma+\sqrt{12} LD^2C^{1/2}}{T^{1/3}}
     \nonumber\\
&  \quad 
+\frac{\sqrt{6} LD^2C^{1/2}}{T^{2/3}}
   + \frac{ D(\sigma^2+G^2)^{1/2}}{T^{1/2}(1-\beta(1-\alpha))}+\frac{LD^2}{2T}
   +\sum_{t=0}^T\left(1-\frac{1}{T^{1/2}}\right)^t \frac{DG}{T}+\sum_{t=0}^T\frac{DG}{T}\left(1-\frac{1}{2T^{2/3}}\right)^{t/2} \nonumber\\
   &\leq \left(1- \frac{1}{T}\right)^T \E{ \frac{1}{n}\sum_{i=1}^n F_i(\bbx^*)- \frac{1}{n}\sum_{i=1}^n F_i(\bar{\bbx}^0)}
+\frac{(1-\alpha)LD^2}{\alpha T} 
+\frac{ LD^2}{T(1-\beta)}    +\frac{\sqrt{6} LD^2C^{1/2}}{T^{2/3}}
     \nonumber\\
&  \quad 
   + \frac{ D(\sigma^2+G^2)^{1/2}}{T^{1/2}(1-\beta(1-\alpha))}+\frac{LD^2}{2T}
   + \frac{DG}{T^{1/2}}+\frac{4DG}{T^{1/3}}   +\frac{\sqrt{2}\sigma+\sqrt{12} LD^2C^{1/2}}{T^{1/3}} ,
\end{align}
%%%
where in the last inequality we use the inequalities $\sum_{t=0}^T\left(1-\frac{1}{2T^{2/3}}\right)^{t/2} \leq  \frac{1}{1-(1-\frac{1}{2T^{2/3}})^{1/2}}\leq 4T^{2/3}$ and $\sum_{t=0}^T\left(1-\frac{1}{T^{1/2}}\right)^t \leq T^{1/2}$. Regrouping the terms and using the inequality $(1-1/T)^T\leq 1/e$ lead to
%%%%
\begin{align} \label{final_thm_discrete_proof_1002}
\E{\frac{1}{n}\sum_{i=1}^n F_i(\bar{\bbx}^{T})}
&\geq (1-e^{-1} ) \frac{1}{n}\sum_{i=1}^n F_i(\bbx^*) 
-\frac{LD^2}{T^{1/2}} 
-\frac{ LD^2}{T(1-\beta)}    -\frac{\sqrt{6} LD^2C^{1/2}}{T^{2/3}}
     \nonumber\\
&  \quad 
   - \frac{ D(\sigma^2+G^2)^{1/2}}{T^{1/2}(1-\beta)}-\frac{LD^2}{2T}
   - \frac{DG}{T^{1/2}}-\frac{4DG}{T^{1/3}}    -\frac{\sqrt{2}\sigma+\sqrt{12} LD^2C^{1/2}}{T^{1/3}}.
\end{align}
%%%%
 Now using the argument in \eqref{final_thm_proof_1500}, we can show that the result in \eqref{final_thm_discrete_proof_1002} implies that for all $j=\ccalN$ it holds
%%%%
\begin{align} \label{final_thm_discrete_proof_10022222}
\E{\frac{1}{n}\sum_{i=1}^n F_i({\bbx_j}^{T})}
&\geq (1-e^{-1} ) \frac{1}{n}\sum_{i=1}^n F_i(\bbx^*) 
-\frac{LD^2}{T^{1/2}} 
-\frac{GD+ LD^2}{T(1-\beta)}    -\frac{\sqrt{6} LD^2C^{1/2}}{T^{2/3}}
     \nonumber\\
&  \qquad 
   - \frac{ D(\sigma^2+G^2)^{1/2}}{T^{1/2}(1-\beta)}-\frac{LD^2}{2T}
   - \frac{DG}{T^{1/2}}-\frac{4DG}{T^{1/3}}  -\frac{\sqrt{2}\sigma+\sqrt{12} LD^2C^{1/2}}{T^{1/3}}.
\end{align}
%%%%
Since $C:=1+\frac{2}{(1-\beta)^2}$ it can be shown that $C^{1/2}=(1+\frac{2}{(1-\beta)^2})^{1/2}\leq 1+\frac{\sqrt{2}}{1-\beta}$. Applying this upper bound into \eqref{final_thm_discrete_proof_10022222} yields the claim in \eqref{eq:main_result_1}.

\bibliography{example_paper,bibliography,bmc,bmc_article_old}
\bibliographystyle{icml2018}

\end{document}